\documentclass{amsart}
\usepackage{amsmath}
\usepackage{epsfig} 
\usepackage{graphics} 
\usepackage{mathrsfs}
\usepackage{amsfonts}
\usepackage{amsfonts,amssymb,amsmath}
\usepackage{epsfig}

\usepackage{amssymb}
\usepackage{tabularx}
\usepackage{enumerate}
\usepackage{graphicx}
\usepackage{texdraw}
\usepackage{dsfont}

\newtheorem{theorem}{Theorem}[section]
\newtheorem{corollary}{Corollary}

\newtheorem{lemma}[theorem]{Lemma}

\theoremstyle{definition}
\newtheorem{definition}[theorem]{Definition}
\newtheorem{remark}{Remark}

\newcommand{\R}{\mathbb{R}}

\newcommand{\wh}{\widehat}

\title[self-similar solutions in a sector]
      {self-similar Solutions in a Sector for a Quasilinear Parabolic Equation$^*$
}

\author[Bendong Lou]{}


\email{blou@tongji.edu.cn}

\thanks{$^*$ This work is supported by NSFC (No. 11271285) and by
the Fundamental Research Funds for the Central Universities (Tongji University).}

\begin{document}
\maketitle

\centerline{\scshape Bendong Lou}
\medskip
{\footnotesize
 \centerline{Department of Mathematics, Tongji University}
   \centerline{ Shanghai 200092, China}
}
\medskip


\medskip

 \centerline{(Communicated by Aim Sciences)}
 \bigskip

\begin{center}
{\it Dedicated to Professor Hiroshi Matano on the occasion of his 60th birthday}
\end{center}
\vskip 5mm

\begin{abstract}
We study a two-point free boundary problem in a sector for a
quasilinear parabolic equation. The boundary
conditions are assumed to be spatially and temporally \lq\lq
self-similar" in a special way. We prove the existence, uniqueness and
asymptotic stability of an expanding solution which is self-similar at discrete
times. We also study the existence and uniqueness of a shrinking solution which is
self-similar at discrete times.
\end{abstract}

\section{Introduction}
Consider the problem
\begin{equation}\label{prob}
\left\{
\begin{array}{l}
 u_t = a(u_x) u_{xx},\quad -\xi_1 (t)<x<\xi_2(t),\ t>0,\\
 u_x(x,t) = -k_1(t,u(x,t)),\ \ u(x,t)=-x\tan\beta\quad \mbox{for }\
x=-\xi_1(t),\ t>0,\\
u_x(x,t) = k_2(t,u(x,t)),\ \ u(x,t)=x\tan\beta\quad \mbox{for }\
x=\xi_2(t),\ t>0,
\end{array}
\right.
\end{equation}
where $a\in C^2(\R)$, $a(\cdot) >0$, $\beta \in (0,\frac{\pi}{2})$
and $k_1,k_2 \in C^2 ([0,\infty)\times [0,\infty), \R)$. In this
problem, $u,\xi_1, \xi_2$ are unknown positive functions to be determined.

The equation in \eqref{prob} includes the heat equation and the curvature flow
equation  as special examples.
In \cite{CGK, CGuo, CGL, GH, GK, Koh}, the authors considered problem \eqref{prob} with
constant $k_i\ (i=1,2)$, that is,
\begin{equation}\label{prob0}
\left\{
\begin{array}{l}
 u_t = a(u_x) u_{xx},\quad -\zeta_1 (t)<x<\zeta_2(t),\ t>0,\\
 u_x(x,t) = -\gamma_1 ,\ \ u(x,t)=-x\tan\beta\quad \mbox{for }\
x=-\zeta_1(t),\ t>0,\\
u_x(x,t) = \gamma_2,\ \ u(x,t)=x\tan\beta\quad \mbox{for }\
x=\zeta_2(t),\ t>0,
\end{array}
\right.
\end{equation}
where $\gamma_1, \gamma_2$ are constants. They proved the existence of solutions of \eqref{prob0}
for some initial data. Moreover, in \cite{CGK, CGL, Koh}, they proved that
when $\gamma_1 +\gamma_2 >0$, any time-global solution $u$ is
expanding (that is, it moves upward to infinity) and it converges asymptotically to
a {\it self-similar solution}: $\sqrt{2t}\; \varphi\left( x/
\sqrt{2t} \right)$. In \cite{CGuo, GH,GK}, they proved that when
$\gamma_1 +\gamma_2 <0$, any solution shrinks to $0$ as
$t\rightarrow T$ for some $T>0$; if $a$ is analytic,
then the rescaled solution $u/\sqrt{2(T-t)}$
converges to a {\it shrinking/backward self-similar solution} with the form: $\psi\big(
x/ \sqrt{2(T-t)} \big)$ as $t\rightarrow T$.

Problem \eqref{prob0} arises in the model of flame propagation in combustion
theory. It also arises in the study of the motion of interface moving with curvature
in which the studied problem is confined in the conical region bounded by two
straight lines and the interface has prescribed touching angles with these two
straight lines (cf. \cite{CGK, CGuo, CGL, GGH, GH, GK, Koh} etc.).
In this paper we will consider a more general problem \eqref{prob}.
In this new problem, the boundary conditions are spatially and temporally
inhomogeneous, which mean that the touching angles between the interface
and the boundaries of the sector domain depend on the spatial and temporal variables.
Clearly, self-similar functions like $\sqrt{2t}\;
\varphi(x/\sqrt{2t})$ or $\sqrt{2(T-t)}\; \psi\big( x/
\sqrt{2(T-t)} \big)$ is no longer a solution of \eqref{prob}. We
have to adopt new concepts for the analogue of self-similar
solutions. Our results in this paper show that problem \eqref{prob}
has an expanding solution which is self-similar at discrete
times if $k_1, k_2$ have some special \lq\lq self-similarity" (see
\eqref{cond:k} below) and if $\min k_1 +\min k_2 >0$. On the other hand,
problem \eqref{prob} has a shrinking solution which is self-similar
at discrete times if $k_1, k_2$ have some special \lq\lq self-similarity" (see
\eqref{cond:k-shrink} or \eqref{cond:k-shrink-equi} below) and if $\max k_1 +\max k_2 <0$.

\begin{definition}\label{def:self1}
Let $(u,\xi_1, \xi_2) = (U, \Xi_1, \Xi_2)$ be a solution of
\eqref{prob} defined for $t\in (0,\infty)$. It is called an {\it expanding self-similar
solution} if
\begin{equation}\label{self1}
b U(x,t) \equiv  U\left( bx, b^2 t  \right) \quad \mbox{for } -
\Xi_1(t) \leqslant x \leqslant \Xi_2(t),\ t>0,
\end{equation}
for some $b>1$, and if
\begin{equation}\label{boundary1}
b \; \Xi_i (t) = \Xi_i \left( b^2 t\right) \quad \mbox{for } t>0
\qquad (i=1,2).
\end{equation}
\end{definition}

From \eqref{self1} we see that, for any $t_0 >0$,
$$
\cdots =  b U(b^{-1} x, b^{-2} t_0 ) = U(x,t_0) = b^{-1} U(bx,b^2
t_0 ) =\cdots.
$$
This means that $U(x,t)$ is similar to $U(x,t_0)$ only at {\it
discrete} times: $t=b^{2m} t_0\ (m\in \mathbb{Z})$. In this sense we
may also say that $(U,\Xi_1, \Xi_2)$ (or, just $U$) is a {\it
discrete expanding self-similar solution} and $\sqrt{2t}\;
\varphi(x/\sqrt{2t})$ is a {\it classical expanding self-similar solution}.

It is easily seen that a necessary condition for the existence of a
discrete expanding self-similar solution is that $k_1$ and $k_2$ are
{\it self-similar} in a special way:
\begin{equation}\label{cond:k}
k_i (t,u)=k_i\left( b^2 t, b u \right) \quad \mbox{for } t,u
\geqslant 0.
\end{equation}
We will give more explanation on this condition near the end of this section.
For simplicity, we also impose another technical conditions on
$k_i$: there exists $\sigma \in (0,\tan\beta)$ such that
\begin{equation}\label{cond:k1}
|k_i(t,u)| \leqslant (\tan\beta)  -\sigma \quad \mbox{for } t, u\geqslant 0 \qquad
(i=1,2).
\end{equation}


\begin{theorem}\label{thm:1}
Assume that $k_1, k_2$ satisfy conditions \eqref{cond:k} and \eqref{cond:k1}.
Assume also that
\begin{equation}\label{cond:expand}
\min k_1 +\min k_2>0
\end{equation}
holds. Then problem \eqref{prob} has a discrete expanding self-similar solution
$(U,\Xi_1, \Xi_2)$.

In addition, if $k_i (t,u)\equiv k_i (u)\ (i=1,2)$, then
\begin{itemize}
\item[\rm (i)] the expanding self-similar solution
is unique and $U_t >0$, $\Xi_{1t} >0$, $\Xi_{2t} >0$ for all
$t>0$;

\item[\rm (ii)] $U$ is asymptotically stable in the sense that
\begin{equation}\label{dH-est}
d_{\mathcal{H}} (\Gamma(t), \gamma(t)) \leqslant C t^{-1/2}\quad \mbox{ as } t\to \infty,
\end{equation}
where $\Gamma(t)$ is the graph of $U$, $\gamma(t)$ is the graph of
any time-global solution $u$ of \eqref{prob} and $d_{\mathcal{H}}$ denotes the Hausdorff
distance.
\end{itemize}
\end{theorem}

\noindent
The existence, uniqueness and asymptotic stability conclusions in this theorem
are proved in subsections 3.6, 3.8 and 3.7, respectively.

Next we consider self-similar solutions which shrink to $0$ in finite
time.

\begin{definition}\label{def:self2}
Given $T>0$. Let $(u,\xi_1, \xi_2) = (\widetilde{U}, \widetilde{\Xi}_1,
\widetilde{\Xi}_2)$ be a solution of \eqref{prob} for $t\in [0,T)$. If
$\|\widetilde{U}(\cdot,t)\|_{L^\infty}\rightarrow 0, \ \widetilde{\Xi}_1
(t)\rightarrow 0, \ \widetilde{\Xi}_2 (t)\rightarrow 0$ as $t\rightarrow T-0$,
\begin{equation}\label{self2}
\widetilde{U}(x,t) \equiv  b\; \widetilde{U}\left( b^{-1}x, b^{-2}t +(1- b^{-2})T
\right) \quad \mbox{for } - \widetilde{\Xi}_1(t) \leqslant x \leqslant
\widetilde{\Xi}_2(t),\ 0\leqslant t <T,
\end{equation}
for some $b>1$, and if
\begin{equation}\label{boundary2}
\widetilde{\Xi}_i (t) = b\; \widetilde{\Xi}_i \left( b^{-2} t +(1-b^{-2})T \right)
\quad \mbox{for } 0\leqslant t<T \qquad (i=1,2),
\end{equation}
then $(\widetilde{U}, \widetilde{\Xi}_1, \widetilde{\Xi}_2)$ (or, just $\widetilde{U}$) is
called a {\it shrinking/backward self-similar solution} of
\eqref{prob} on time interval $[0,T)$.
\end{definition}

Since $\widetilde{U}(x,t)$ is similar to $\widetilde{U}(x,t_0)$ only at discrete
times: $t=b^{-2m} t_0 +(1-b^{-2m})T$  ($m\in \mathbb{Z}$ and
$2m\geqslant \log(\frac{T-t_0}{T})/\log b$), we may also say that
$(\widetilde{U}, \widetilde{\Xi}_1, \widetilde{\Xi}_2)$  is a {\it discrete shrinking
self-similar solution} on $[0,T)$.

A necessary condition for the existence of such a solution is that
\begin{equation}\label{cond:k-shrink}
k_i (t,u) = k_i\left(b^{-2} t+(1-b^{-2})T, \; b^{-1} u \right) \quad
\mbox{for } 0\leqslant t <T,\ u \geqslant 0.
\end{equation}
Replacing $t$ by $T-t'$ then we see that \eqref{self2}, \eqref{boundary2}
and \eqref{cond:k-shrink} are equivalent to
\begin{equation}\label{self2-equi}
\widetilde{U}(x,T-t') \equiv  b \; \widetilde{U}\left( b^{-1}x, T- b^{-2}t'\right)
\end{equation}
for $- \widetilde{\Xi}_1(T-t') \leqslant x \leqslant \widetilde{\Xi}_2(T-t')$,
$0< t' \leqslant T$,
\begin{equation}\label{boundary2-equi}
\widetilde{\Xi}_i (T-t') = b\; \widetilde{\Xi}_i \left( T- b^{-2} t' \right) \quad
\mbox{for } 0<t' \leqslant T \qquad (i=1,2),
\end{equation}
and
\begin{equation}\label{cond:k-shrink-equi}
k_i (T-t' , u) = k_i (T- b^{-2} t' , b^{-1} u) \quad \mbox{for }
0<t' \leqslant T,\ u\geqslant 0 \qquad (i=1,2),
\end{equation}
respectively.

\begin{theorem}\label{thm:2}
Given $T>0$, assume that $k_1, k_2$ satisfy condition \eqref{cond:k-shrink}
or \eqref{cond:k-shrink-equi}. Assume also that \eqref{cond:k1} and
\begin{equation}\label{cond:shrink}
\max k_1 +\max k_2 <0
\end{equation}
hold. Then problem \eqref{prob} has a discrete shrinking self-similar
solution $(\widetilde{U}, \widetilde{\Xi}_1, \widetilde{\Xi}_2)$ on $[0,T)$.

In addition, if $k_i(t,u) \equiv k_i (u)\ (i=1,2)$, then the discrete shrinking
self-similar solution is unique and $\widetilde{U}_t <0,\
\widetilde{\Xi}_{1t} <0$, $\widetilde{\Xi}_{2t} <0$ for $t\in [0,T)$.
\end{theorem}

The uniqueness for shrinking self-similar solutions is not necessary to be true,
even for the special problem \eqref{prob0} (cf. \cite{GH, GK})).
But the above theorem shows that it can be unique under certain assumptions.

Definition \ref{def:self2} and Theorem \ref{thm:2} deal with shrinking
solutions on finite time interval $[0,T)$. If we take a time shift, these
solutions can be regarded as solutions defined on $[-T,0)$. More precisely,
let $(\widetilde{U},\widetilde{\Xi}_1, \widetilde{\Xi}_2)$ be a shrinking self-similar
solution of \eqref{prob} on $[0,T)$. Then
$$
\wh{U}(x,t;T) := \widetilde{U}(x,T+t)\quad \mbox{for } -\wh{\Xi}_1
(t;T)\leqslant x \leqslant \wh{\Xi}_2(t;T),\ t\in [-T, 0),
$$
and
$$
\wh{\Xi}_i (t;T) := \widetilde{\Xi}_i (T+t)\quad \mbox{for }  t\in [-T,
0)\quad (i=1,2)
$$
satisfy
\begin{equation}\label{self3}
\wh{U}(x,t) \equiv  b \; \wh{U} ( b^{-1}x, b^{-2}t ) \quad \mbox{for }
- \wh{\Xi}_1(t) \leqslant x \leqslant \wh{\Xi}_2 (t),\ -T\leqslant
t<0,
\end{equation}
and
\begin{equation}\label{boundary3}
\wh{\Xi}_i (t) = b\; \wh{\Xi}_i ( b^{-2} t ) \quad \mbox{for } -T
\leqslant t<0 \qquad (i=1,2).
\end{equation}
So $(\wh{U},\wh{\Xi}_1, \wh{\Xi}_2)$ (which is defined on $[-T,0)$
and shrinks to $0$ as $t\rightarrow 0-0$) is a self-similar solution
of
\begin{equation}\label{prob-neg}
\left\{
\begin{array}{l}
u_t = a(u_x) u_{xx},\quad -\xi_1 (t)<x<\xi_2(t),\ t<0,\\
u_x(x,t) = -k_1(T+t, u(x,t)),\ \ u(x,t)=-x\tan\beta\quad \mbox{for
}\ x=-\xi_1(t),\ t<0,\\
u_x(x,t) = k_2(T+t, u(x,t)),\ \ u(x,t)=x\tan\beta\quad \mbox{for }\
x=\xi_2(t),\ t<0.
\end{array}
\right.
\end{equation}

We now consider shrinking self-similar solutions defined in
$(-\infty,0)$.

\begin{definition}\label{def:self3}
Assume that $k_i (t,u)\equiv k_i (u)\ (i=1,2)$. Let $(\wh{U}, \wh{\Xi}_1,
\wh{\Xi}_2)$ be a solution of \eqref{prob-neg} defined for $t\in (-\infty,
0)$. If it satisfies \eqref{self3} and \eqref{boundary3} for some
$b>1$ and $t\in (-\infty,0)$, then it is called a {\it discrete shrinking/backward
self-similar solution} of \eqref{prob-neg} in
$(-\infty,0)$.
\end{definition}

\begin{theorem}\label{thm:3}
Assume that $k_1\equiv k_1(u)$ and $k_2\equiv k_2(u)$ satisfy \eqref{cond:k1},
\eqref{cond:shrink} and $k_i (u)=k_i (bu)$ for all $u\geqslant 0$ and some $b>1$.
Then problem \eqref{prob-neg} has a unique discrete shrinking self-similar solution
$(\wh{U}, \wh{\Xi}_1, \wh{\Xi}_2)$ in $(-\infty, 0)$, and
$\wh{U}_t <0,\ \wh{\Xi}_{1t} <0, \ \wh{\Xi}_{2t}<0$ for $t\in (-\infty, 0)$.
\end{theorem}

Our theorems extend the results about classical self-similar
solutions in \cite{CGK,  CGuo, CGL, GH, GK, Koh} to problem \eqref{prob} with
nonlinear boundary conditions. Our approach is essentially different from theirs
though we will use their classical self-similar solutions as lower
and upper solutions to give the growth bound for the solution of
\eqref{prob}. We will convert probe \eqref{prob} by changing
variables to a new problem in a fixed domain. Then we use a convergence result in
\cite{BPS} to show that the $\omega$-limit of the unknown in the new problem
is a periodic solution, which corresponds to a discrete self-similar solution of
\eqref{prob}.

Our boundary conditions are given by functions $k_1$ and $k_2$ which are self-similar
as in \eqref{cond:k}. We now give some examples and/or backgrounds on such kind of
self-similarity. First, some reactions in chemistry occur in a media with obstacles (cf. \cite{ten}).
When the obstacles arrange in a regular way, it is possible to be studied
from a mathematical point of view.
For example, if we consider a Belousov-Zhabotinsky (BZ) reaction in a media with obstacles
arranging in columns, then the interface propagation in the BZ experiment can be studied
through a curvature flow in a band domain with undulating boundaries (cf. \cite{LMN, MNL}).
Similarly, if we consider the BZ reaction in a media with obstacles arranging in radial
rays with center at origin $O$,  and if the ratios of the sizes of adjacent obstacles
are constant, then the interface propagation can be studied through a
curvature flow in a sector with undulating boundaries, which is essentially a similar
problem as our \eqref{prob}.
Another example is the following. In geology, Liesegang rings are colored bands of
cement observed in sedimentary rocks, which are often referred to as great examples
of geochemical self-organization (cf. \cite{stow}). Generally, the Liesegang rings are
arranged in a regular self-similar way: the ratios of the widths of adjacent annuluses
are constant (cf. \cite{HHMO, KR, stow}). If we cut off a sector with apex at the
center of the rings and consider the interface propagation in this notch, then the
problem is reduced to one like \eqref{prob}.


In Section 2 we give some preliminaries, including the selection of
the initial data, the local existence result and comparison principles.
In Section 3 we consider the expanding case and prove
Theorem \ref{thm:1}. In Section 4 we consider the shrinking case and prove
Theorems \ref{thm:2} and \ref{thm:3}.


\section{Preliminaries}
We use notation $S := \{(x,y)\mid y> |x|\tan \beta,\ x\in \R\}$, and
use $\partial_1 S$ and $\partial_2 S$ to denote the left and right
boundaries of $S$, respectively.
For any $t_0 >0$, let $(u(x,t), \xi_1 (t),\xi_2(t))$ be a classical
solution of \eqref{prob} on the time interval $[0,t_0]$ with some initial data.
Then we write
$$
Q_{t_0} := \{(x,t)\; |\; -\xi_1 (t) < x < \xi_2 (t)\ {\rm and}\ 0< t
\leqslant t_0\}.
$$

\subsection{Initial data}
We will consider the problem \eqref{prob} with initial data
\begin{equation}\label{ini}
u(x,0)=u_0(x),\quad -\xi_{01}\leqslant x\leqslant \xi_{02},
\end{equation}
where $u_0(x)>0$, $\xi_{01}>0$ and $\xi_{02}>0$ satisfy
\begin{equation}\label{ini-boundary}
u_0 (-\xi_{01}) = \xi_{01} \tan\beta , \quad
u_0 (\xi_{02}) = \xi_{02} \tan\beta
\end{equation}
and the compatibility conditions:
\begin{equation}\label{compatibility}
(u_0)_x (-\xi_{01}) = -k_1 (0,u_0 (-\xi_{01})), \quad (u_0)_x
(\xi_{02}) = k_2 (0,u_0 (\xi_{02})).
\end{equation}

Since our main purpose in this paper is to construct self-similar
solutions, we will not focus on general solutions of \eqref{prob} and \eqref{ini} for
general $u_0$ as it was done in \cite{CGK, GH, Koh}, but choose $u_0 \in C^{2+\mu} ([-\xi_{01},\xi_{02}])$
for some $\mu\in (0,1)$, and only consider classical solution $u$ of \eqref{prob} and \eqref{ini}
in $C^{2+\mu, 1+\mu/2}(\overline{Q}_{t_0})$ for $t_0 >0$. Moreover, we require that $u_0$ satisfies
\begin{equation}\label{ini-gradient}
|(u_0)_x (x)| \leqslant \tan \beta -\sigma
\quad \mbox{for } -\xi_{01}\leqslant x \leqslant \xi_{02},
\end{equation}
where $\sigma$ is as in \eqref{cond:k1}. This inequality does not conflict with
the compatibility conditions by \eqref{cond:k1}.

In summary, in this paper we choose
initial data from the following set of admissible functions:
\begin{equation}\label{admissiable}
C^{2+\mu}_{\rm ad} := \left\{ u_0 \left|
\begin{array}{l}
 u_0 \in C^{2+\mu}([-\xi_{01},\xi_{02}]) \mbox{ for some } \mu\in (0,1), \ \mbox{where } \\
 \xi_{01},\; \xi_{02} >0 \mbox{ and } u_0(\cdot )>0  \mbox{ satisfy } \eqref{ini-boundary},\
 \eqref{compatibility} \mbox{ and }  \eqref{ini-gradient}
 \end{array}
\right. \right\}.
\end{equation}

\subsection{Gradient bound of $u$}

\begin{lemma}\label{lem:gradient bound}
    Let $u_0 \in C^{2+\mu}_{\rm ad}$ for some $\mu\in (0,1)$,
    $u(x,t)\in C^{2+\mu,1+\mu/2} (\overline{Q}_{t_0})$ be a
    solution of  \eqref{prob} and \eqref{ini} on $[0,t_0]$. Then
\begin{equation}\label{gradient bound}
 |u_x (x,t)| \leqslant \tan\beta - \sigma \quad \mbox{for }
 (x,t)\in  \overline{Q}_{t_0}.
\end{equation}
\end{lemma}

\begin{proof}
By \eqref{cond:k1} we have
$$
u_x (\xi_2 (t),t) = k_2 (t,u (\xi_2 (t),t) ) \leqslant \tan\beta
-\sigma,
$$
and
$$
u_x (-\xi_1 (t),t) = -k_1 (t,u (-\xi_1 (t),t) ) \leqslant \tan\beta -\sigma.
$$
Combining these inequalities with \eqref{ini-gradient} we obtain
$u_x \leqslant \tan\beta -\sigma$ by maximum principle. $u_x
\geqslant -\tan\beta +\sigma$ is proved similarly.
\end{proof}

\begin{corollary}
       Let $u_0$ and $u$ be as in the previous lemma. Then, for
    any $\theta\in [-\theta_0, \theta_0]$ with $\theta_0 := \frac{\pi}{2}
    -\beta$, there holds
\begin{equation}\label{regular1}
 \sigma\cos\beta \leqslant    (1 \pm u_x (x,t) \tan\theta) \cos\theta \leqslant
 2- \sigma\cot\beta \quad \mbox{for }  (x,t)\in  \overline{Q}_{t_0}.
\end{equation}
\end{corollary}

\subsection{Change of variables}

To study the local and global existence of solutions of the initial
boundary value problem \eqref{prob} and \eqref{ini}, it is convenient to
introduce new coordinates that convert the sector domain $S$ into a flat
cylinder.  More precisely, we will make a change of variables
$(x,y,t)\mapsto (\theta,\rho,s)$, which gives a diffeomorphism
$(\overline{S}\backslash \{0\})\times [0,t_\infty)\to
\overline{D}\times [s_0,s_\infty)$, where
$$
 D:=\{(\theta,\rho)\in\R^2 \mid -\theta_0<\theta<\theta_0,\
 -\infty<\rho<\infty\}
$$
with $\theta_0 := \frac{\pi}{2} -\beta$. The functions
$\theta=\theta(x,y,t)$, $\rho= \rho(x,y,t)$ and $s=s(t)$ are to be
specified below. With these new coordinates, the function $y=u(x,t)$
is expressed as $\rho=\omega (\theta,s)$, where the new unknown
$\omega(\theta,s)$ is determined by the relation
\begin{equation}\label{u-to-v}
 \rho\left(x,u(x,t),t\right)=\omega\left(\theta(x,u(x,t),t),s(t)\right).
\end{equation}
The function $\omega(\theta,s)$ is well-defined
provided that the map $t\mapsto s(t)$ is strictly monotone for
$t\in [0, t_\infty )$ and $x\mapsto \theta(x,u(x,t),t)$ is strictly
monotone for each fixed $t\in [0, t_\infty)$. We will see later that these
monotonicity conditions always hold for the class of solutions that
we consider. Indeed we will  prove
\begin{equation}\label{theta-u}
 \frac{\partial}{\partial t} s(t) > 0,\quad
 \frac{\partial }{\partial x}\theta \left(x,u(x,t),t\right)
 =\theta_x+\theta_y u_x >0.
\end{equation}

Once $\omega(\theta,s)$ is defined, then substituting it into the
relation $y=u(x,t)$ yields
\begin{equation}\label{v-to-u}
 Y\left(\theta,\omega(\theta,s),s\right)=u\left(X(\theta,\omega(\theta,s),s),T(s)\right),
\end{equation}
where the map $(\theta,\rho,s)\mapsto\left( X(\theta,\rho,s),
Y(\theta,\rho,s),T(s)\right): \overline{D}\times [s_0,s_\infty) \to
(\overline{S}\backslash \{0\})\times [0,t_\infty)$ is the inverse
map of $(x,y,t) \mapsto \left( \theta(x,y,t), \rho(x,y,t), s(t)
\right)$. The expression \eqref{v-to-u} gives a formula for
recovering the original solution $u(x,t)$ from $\omega(\theta,s)$.
In order for $u$ to be smoothly dependent on $\omega$, we need the
map $\theta\mapsto X(\theta,\omega(\theta,s),s)$ to be one-to-one
for each fixed $s$ and that $s\mapsto T(s)$ is strictly monotone for
$s\in [s_0, s_\infty)$. Indeed we will prove
\begin{equation}\label{X-v}
\frac{\partial }{\partial s} T(s) > 0,\quad
 \frac{\partial }{\partial \theta}X\left(\theta, \omega(\theta,s),s\right)
 =X_\theta+X_\rho \omega_\theta >0.
\end{equation}

\subsection{Local existence}\label{local}

To get the local existence we make the following change of
variables.
\begin{equation}\label{new-variable-1}
\left\{
\begin{array}{ll}
\displaystyle \theta = \arctan \frac{x}{y}, & (x,y)\in
\overline{S}\backslash \{0\}, \\
\displaystyle \rho= \frac{1}{2} \log (x^2 +y^2), &
(x,y)\in \overline{S} \backslash \{0\},\\
\displaystyle s= t,& t\geqslant 0.
\end{array}
\right.
\end{equation}
The inverse map is
\begin{equation}\label{inverse-map}
\left\{
\begin{array}{ll}
x= e^\rho \sin \theta,& (\theta, \rho)\in \overline{D},\\
y= e^\rho \cos \theta,& (\theta, \rho)\in \overline{D},\\
t = s, & s\in [0,\infty).
\end{array}
\right.
\end{equation}
Clearly, $\theta =\theta_0$ and $\theta = -\theta_0$ correspond to
$\partial_2 S$ and $\partial_1 S$, respectively.

Let $u(x,t) >0$ be a classical solution of \eqref{prob} and \eqref{ini}
for $t\geqslant 0$, then
\begin{equation}\label{omega-def}
\rho(x, u(x,t),t) = \omega (\theta (x,u(x,t),t),s)\quad
\Leftrightarrow \quad e^{\omega (\theta,s)} \cos \theta = u \Big(
e^{\omega(\theta,s)}\sin \theta, t \Big).
\end{equation}
defines a new unknown $\rho=\omega (\theta,s)$ for $s\geqslant 0$.
This function is well-defined since
$$
\frac{\partial}{\partial x} \theta(x,u(x,t),t)=
\frac{\partial}{\partial x} \left( \arctan \frac{x}{u(x,t)} \right)
=\frac{1- u_x \tan\theta}{u\cdot (1+\tan^2 \theta)} >0
$$
by \eqref{regular1}.

Differentiating the expression $e^{\omega (\theta,s)} \cos \theta =
u ( e^{\omega(\theta,s)}\sin \theta, t)$ twice by $\theta$ and once
by $t$ we obtain
$$
u_x = \frac{\omega_\theta \cos \theta - \sin\theta}{\cos \theta
+\omega_\theta \sin\theta }, \quad  u_{xx} =
\frac{\omega_{\theta\theta} - \omega^2_\theta -1} { e^\omega (\cos
\theta + \omega_\theta \sin\theta)^3 }, \quad u_t = \frac{e^\omega
\omega_s }{\cos \theta + \omega_\theta \sin\theta}.
$$
Therefore, problem \eqref{prob} with \eqref{ini} is converted into the
following problem
\begin{equation}\label{prob-omega}
\left\{
\begin{array}{ll}
\displaystyle \omega_s = a\left( \frac{\omega_\theta \cos \theta -
\sin \theta} {\cos\theta +\omega_\theta \sin \theta}\right)
\frac{\omega_{\theta\theta} - \omega^2_\theta -1}{e^{2\omega}
(\cos\theta +\omega_\theta \sin \theta)^2} , & \theta\in (-\theta_0 ,\theta_0),
\ s\in (0,\infty),\\
\omega_\theta (-\theta_0,s) = - h^0_1 (s,\omega(-\theta_0 ,s)), & s\in [0,\infty),\\
\omega_\theta (\theta_0,s) = h^0_2 (s,\omega(\theta_0 ,s)), & s\in
[0,\infty),\\
\omega(\theta, 0)= \tilde{\omega} (\theta), & \theta\in [-\theta_0,
\theta_0].
\end{array}
\right.
\end{equation}
where $\tilde{\omega}$ is defined by \eqref{omega-def} at $t=s=0$ and
\begin{equation}\label{def:h-rho}
h^0_i (s,\omega)=\frac{\sin\theta_0 + k_i (s, e^\omega \cos
\theta_0 ) \cos\theta_0}{\cos\theta_0- k_i (s, e^\omega \cos
\theta_0 ) \sin\theta_0} \quad (i=1,2).
\end{equation}

Estimate \eqref{gradient bound} implies that
\begin{equation}\label{est omega}
\sigma - \tan\beta \leqslant u_x = \frac{\omega_\theta \cos\theta-
\sin \theta}{\cos\theta + \omega_\theta \sin \theta} \leqslant
\tan\beta -\sigma \quad \mbox{for } \theta\in [-\theta_0,
\theta_0],\ s\geqslant 0.
\end{equation}
Thus $\cos\theta + \omega_\theta \sin\theta > 0$ since it is positive at $\theta =0$
and it can not be zero by \eqref{est omega}. Considering the second
inequality in \eqref{est omega} we have
\begin{equation}\label{est omega-theta}
\omega_\theta [\cos\theta - \sin\theta (\tan\beta -\sigma)]
\leqslant \sin \theta + (\tan\beta -\sigma)\cos\theta.
\end{equation}
Note that, for $\theta \in [-\theta_0, \theta_0]$ ($\theta_0 =
\frac{\pi}{2}-\beta$), we have
$$
\cos\theta - \sin\theta (\tan\beta -\sigma) \geqslant \cos\theta_0 [
1- \tan\theta_0 (\tan\beta -\sigma)] \geqslant \sigma \cos\beta,
$$
$$
\cos\theta - \sin\theta (\tan\beta -\sigma) \leqslant \cos\theta [ 1
+ \tan\theta_0 (\tan\beta -\sigma)] \leqslant 2-\sigma \cot \beta.
$$
So
$$
\omega_\theta \leqslant \frac{\sin \theta + (\tan\beta
-\sigma)\cos\theta}{\cos\theta -\sin\theta (\tan\beta -\sigma)}
\leqslant \Omega_1 := \frac{1 + \tan\beta -\sigma}{\sigma \cos \beta
}.
$$
Using the first inequality in this formula we have, for $\theta \leqslant 0$,
$$
\cos \theta + \omega_\theta \sin \theta \geqslant
\frac{1}{\cos\theta -\sin\theta(\tan\beta -\sigma)} \geqslant
\varepsilon_1 := \frac{1}{2-\sigma \cot \beta}.
$$

Similarly, considering the first inequality in \eqref{est omega} we
have
$$
\omega_\theta \geqslant \frac{\sin \theta - (\tan\beta -\sigma)
\cos\theta}{\cos \theta + \sin\theta (\tan\beta -\sigma)} \geqslant
- \Omega_1 ,
$$
and $\cos \theta +\omega_\theta \sin \theta \geqslant \varepsilon_1$ for
$\theta\geqslant 0$.

Summarizing the above results we have
\begin{equation}\label{gradient omega}
|\omega_\theta (\theta, s)| \leqslant \Omega_1 \quad \mbox{and}
\quad \cos \theta + \omega_\theta \sin \theta \geqslant
\varepsilon_1 >0
\end{equation}
for $\theta \in [-\theta_0, \theta_0]$ and $s\geqslant 0$.

By the standard theory for parabolic equations, we see that
\eqref{prob-omega} has a classical solution on time interval $s\in
[0,2\tau]$ for positive $\tau = \tau(k_1, k_2, \mu, \tilde{\omega})$.

The second inequality in \eqref{gradient omega} implies that, once
the solution $\omega$ of \eqref{prob-omega} is obtained then we can
recover it to the original solution $u$ of \eqref{prob}. In fact,
$$
\frac{\partial }{\partial \theta} X(\theta, \omega(\theta,s),s) =
\frac{\partial }{\partial \theta} e^{\omega (\theta,s)} \sin \theta
= e^{\omega (\theta,s)} (\cos\theta +\omega_\theta \sin\theta)
>0.
$$
Consequently, we have the following local existence result.

\begin{lemma}\label{lem:local}
Problem \eqref{prob} with initial data $u_0(x)\in C^{2+\mu}_{\rm ad}\
(\mu \in (0,1))$ has a classical solution
$u$ on time interval $[0,2\tau]$, where $\tau$ depends only on $k_1,
k_2, \mu$ and $u_0$.
\end{lemma}

\subsection{Comparison principle}

Let $v_1(x), v_2 (x)$ be two functions whose graphs lie in $\overline{S}$ and meet
the two boundaries of $S$. Hereafter, when we write
$$
v_1 \leqslant v_2  \qquad (\mbox{resp. } v_1 \ll v_2\; ),
$$
we mean that $v_1 (x)\leqslant v_2 (x)$ (resp. $v_1(x) <v_2(x)$) for all $x$
with $(x,v_i(x))\in \overline{S}\ (i=1,2)$; when we write
$$
v_1  \preceq v_2
$$
we mean that $v_1 (x) \leqslant v_2 (x)$ and the \lq\lq equality" holds at
some $x$.

Assume further that $|v_{1x}|, |v_{2x}| <\tan\beta$. Then for each $x$ with
$(x,v_1(x))\in \overline{S}\backslash \{O\}$,
there exists a unique $Z(x)$ such that
$$
x\cdot v_2 (Z(x)) = Z(x) \cdot v_1(x),
$$
that is, $(x,v_1(x))$ and $(Z(x), v_2 (Z(x)))$ lie on the same line passing the origin.
By a simple geometric observation we have
\begin{equation}\label{op-1}
v_1 \leqslant  v_2
\ \Leftrightarrow \ x^2 +v^2_1(x) \leqslant Z^2(x) + v_2^2 (Z(x))\mbox{ for }
x \mbox{ with } (x,v_1(x))\in \overline{S}\backslash \{O\}
\end{equation}
and
\begin{equation}\label{op-2}
v_1 \ll   v_2 \ \Leftrightarrow \ x^2 +v^2_1(x) < Z^2(x) + v_2^2 (Z(x))\mbox{ for }
x \mbox{ with } (x,v_1(x))\in \overline{S}\backslash \{O\}.
\end{equation}

For some $t_0 >0$, let $u_1 (x,t)\in C^{2,1}\big(\overline{Q_{t_0}^{(1)}}\big)$ and
$u_2 (x,t)\in C^{2,1}\big(\overline{Q_{t_0}^{(2)}}\big)$ be two positive
functions, where, for $i=1,2$,
$$
Q_{t_0}^{(i)}:= \{(x,t)\mid -\xi_1^{(i)}(t) <x< \xi_2^{(i)}(t),\ 0<t\leqslant t_0\},
$$
$$
u_i (-\xi_1^{(i)}(t), t) =\xi_1^{(i)}(t) \cdot \tan\beta,\ \ u_i (\xi_2^{(i)}(t), t)
= \xi_2^{(i)}(t) \cdot \tan\beta,\quad 0<t\leqslant t_0,
$$
and $|(u_i)_x (x,t)| < \tan\beta$.

\begin{definition}
Let $t_0 >0$  and $u_i \in C^{2,1} \big(\overline{Q_{t_0}^{(i)}}\big) \ (i=1,2)$ be positive functions
as above. Then $u_1$ is called a lower solution of \eqref{prob} on $[0,t_0]$ if
 \begin{equation}\label{u1}
\left\{
 \begin{array}{l}
  u_{1t} \leqslant a(u_{1x})  u_{1xx} \ \ \ \ \ \ \ {\rm for }\ -\xi_1^{(1)}(t) <x <
  \xi_2^{(1)}(t),\ 0\leqslant t\leqslant t_0,\\
      u_{1x} (x, t) \geqslant  -k_1 (t, u_1 (x,t))
       \ \ \ \ \ \  {\rm for}\ x= - \xi_1^{(1)} (t),\ 0\leqslant t\leqslant t_0,\\
      \displaystyle
      u_{1x} (x, t) \leqslant  k_2 (t, u_1 (x,t)) \ \ \ \ \ \
      {\rm for}\ x=\xi_2^{(1)} (t),\ 0\leqslant t\leqslant t_0.
     \end{array}
     \right.
     \end{equation}
Similarly, $u_2$ is called an upper solution of \eqref{prob} if the opposite
inequalities hold.
\end{definition}

The following comparison principle holds.

\begin{lemma}\label{comparison}
      Let $t_0 >0$. Assume that $u_1 (x,t)$ and $u_2 (x,t)$ are lower solution and
      upper solution of \eqref{prob} on $[0, t_0]$, respectively.
      If $u_1 (\cdot,0) \leqslant u_2 (\cdot,0)$, then $u_1 (\cdot, t) \leqslant u_2 (\cdot,t)$
      for $0\leqslant t \leqslant t_0$.
     If $u_1 (\cdot,0) \leqslant u_2 (\cdot,0)$ and $u_1 (x,0)\not\equiv u_2(x,0)$,
     then $u_1 (\cdot, t) \ll  u_2 (\cdot,t)$ for $0< t \leqslant t_0$.
\end{lemma}

\begin{proof}
We change variables by \eqref{new-variable-1} and \eqref{inverse-map}, that is,
using
$$
e^{\rho_i } \cos \theta = u_i ( e^{\rho_i }\sin \theta, s ),
$$
we define implicit functions $\rho_i = \omega_i (\theta,s)\ (i=1,2)$.
Since $u_1(x,t)$ is a lower solution of \eqref{prob}, it is easily seen that
$\omega_1 (\theta,s)$ is a lower solution of \eqref{prob-omega}:
$$
\left\{
\begin{array}{ll}
\displaystyle \omega_{1s} \leqslant a\left( \frac{\omega_{1\theta} \cos \theta -
\sin \theta} {\cos\theta +\omega_{1\theta} \sin \theta}\right)
\frac{\omega_{1\theta\theta} - \omega^2_{1\theta} -1}{e^{2\omega_1}
(\cos\theta +\omega_{1\theta} \sin \theta)^2} , & \theta\in (-\theta_0 ,\theta_0),
\ s\in (0,t_0],\\
\omega_{1\theta} (-\theta_0,s) \geqslant - h^0_1 (s,\omega_1(-\theta_0 ,s)), & s\in [0,t_0],\\
\omega_{1\theta} (\theta_0,s) \leqslant h^0_2 (s,\omega_1(\theta_0 ,s)), & s\in [0,t_0].
\end{array}
\right.
$$
Similarly, $\omega_2(\theta,s)$ is an upper solution of \eqref{prob-omega}.
By \eqref{op-1}, $u_1(\cdot, t)\leqslant u_2(\cdot, t)$ is equivalent to
$\omega_1 (\theta,s) \leqslant \omega_2(\theta,s)$.  The latter follows from
the comparison principle for \eqref{prob-omega}, which is a problem in a fixed domain.
Similarly, the conclusion $u_1(\cdot,t)\ll u_2(\cdot,t)$ can be proved by using \eqref{op-2}.
\end{proof}


\section{Expanding self-similar solutions}
\par
In this section we always assume that \eqref{cond:expand} holds.

\subsection{Classical expanding self-similar solutions}
We will use classical self-similar solutions of \eqref{prob0} as upper and lower solutions of
\eqref{prob} to give the growth bound for the solution $u$ of
\eqref{prob} and \eqref{ini}.
For any $\gamma_1, \gamma_2 \in \R$, consider the problem
\begin{equation}\label{eq-phi}
\left\{
\begin{array}{l}
a(\varphi'(z))\varphi''(z) = \varphi (z) -z\varphi'(z),\ \ \ z\in \R,\\
\varphi'(-p_1) = -\gamma_1, \ \varphi(-p_1) = p_1 \tan\beta,\\
\varphi'(p_2) = \gamma_2, \ \varphi(p_2) = p_2\tan\beta.
\end{array}
\right.
\end{equation}
In \cite{CGK, CGL, Koh}, the authors obtained the following result.

\begin{lemma}
For any given $\gamma_1, \gamma_2$ with $\gamma_1 + \gamma_2 >0$,
there exists a unique pair $p_1, p_2 >0$ such that problem
\eqref{eq-phi} has a solution $\varphi(z;\gamma_1, \gamma_2)$, which
is positive on $[-p_1, p_2]$.
\end{lemma}

It is easily seen that the function
$$
\sqrt{2t}\; \varphi\left( \frac{x}{ \sqrt{2t} }; \gamma_1,
\gamma_2\right) \quad \mbox{for } - \zeta_1(t) < x < \zeta_2(t),\
t>0,
$$
with $\zeta_i (t) = p_i \sqrt{2t}\ (i=1,2)$ is an expanding self-similar solution of \eqref{prob0}.
Set
\begin{equation}\label{cond:k2}
k^0_i := \min k_i(t,u) \quad \mbox{and} \quad K^0_i := \max k_i (t,u)
\qquad (i=1,2),
\end{equation}
and define
$$
\varphi^- (z):= \varphi(z; k^0_1, k^0_2), \quad
\varphi^+ (z) := \varphi(z; K^0_1, K^0_2).
$$
Then $\sqrt{2t}\; \varphi^- (x/\sqrt{2t})$ and $\sqrt{2t}\; \varphi^+ (x/\sqrt{2t})$
(both are expanding self-similar solutions of \eqref{prob0}) are
lower and upper solutions of \eqref{prob}, respectively.

Since the initial data $u_0 >0$, there exist $t^+, t^- >0$ such
that
\begin{equation}\label{bound of u0}
\sqrt{2t^-}\; \varphi^-  \left({\frac {\cdot } {\sqrt{2t^-}}}\right)
\leqslant u_0(\cdot ) \leqslant \sqrt{2t^+}\; \varphi^+  \left({\frac {\cdot }
{\sqrt{2t^+}}}\right).
\end{equation}
The comparison principle implies that
\begin{equation}\label{bound of u}
\sqrt{2(t+t^-)}\; \varphi^-  \left({\frac {\cdot } {\sqrt{2(t+
t^-)}}}\right) \leqslant u(\cdot ,t) \leqslant \sqrt{2(t+t^+)}\;
\varphi^+ \left({\frac {\cdot } {\sqrt{2(t+t^+)}}}\right).
\end{equation}

\subsection{Changes of variables}

In subsection \ref{local} we gave a local existence result. One
difficulty for deriving the global existence is the lack of the
growth bound for $\omega$. To give the global existence we adopt another
change of variables.

Let $\tau$ be the constant in Lemma \ref{lem:local} and let $t^- >0$
be as in \eqref{bound of u0}. For any $n\in \mathbb{N}$
satisfying
\begin{equation}\label{choose of n}
n>\frac{1}{t^-}\quad \mbox{and } \quad n> \frac{1}{\tau},
\end{equation}
we introduce new variables by
\begin{equation}\label{new-variable-v}
\left\{
\begin{array}{ll}
\displaystyle \theta = \arctan \frac{x}{y}, & (x,y)\in
\overline{S}\backslash
\{0\}, \\
\displaystyle \rho= \frac{1}{2} \log \frac{n(x^2 +y^2)}{nt+ 1}, &
(x,y)\in \overline{S} \backslash \{0\},\ t\geqslant 0,\\
\displaystyle s= \frac{1}{2} \log \Big(t +\frac{1}{n} \Big),&
t\geqslant 0.
\end{array}
\right.
\end{equation}
The inverse map is
\begin{equation}\label{inverse-map-v}
\left\{
\begin{array}{ll}
x= e^s e^\rho \sin \theta,& (\theta, \rho)\in \overline{D},\ s\in [-\frac{1}{2}\log n,\infty),\\
y=e^s e^\rho \cos \theta,& (\theta, \rho)\in \overline{D},\ s\in [-\frac{1}{2}\log n,\infty),\\
t = e^{2s}-\frac{1}{n}, & s\in [-\frac{1}{2}\log n,\infty).
\end{array}
\right.
\end{equation}
Clearly, $\theta =\theta_0$ and $\theta = -\theta_0$ correspond to
$\partial_2 S$ and $\partial_1 S$, respectively.

Let $u(x,t) >0$ be a solution of \eqref{prob} and \eqref{ini} for
$t\geqslant 0$, then a similar discussion as in subsection
\ref{local} shows that
\begin{equation}\label{v-def-1}
\rho(x, u(x,t),t) = v (\theta (x,u(x,t),t),s(t)) \
\Leftrightarrow \ e^s e^{v (\theta,s)} \cos \theta = u \Big( e^s
e^{v(\theta,s)}\sin \theta, e^{2s}-\frac{1}{n} \Big)
\end{equation}
defines a new unknown $\rho=v (\theta,s)$ for $s\in
[-\frac{1}{2}\log n, \infty)$. Differentiating the second equality
twice by $\theta$ and once by $s$ we obtain
$$
u_x = \frac{v_\theta \cos \theta - \sin\theta}{\cos \theta +v_\theta
\sin\theta },\quad u_{xx} = \frac{v_{\theta\theta} -v^2_\theta -1}
{e^s e^v (\cos \theta +v_\theta \sin\theta)^3 }, \quad u_t =
\frac{e^v (1+v_s)}{2e^s (\cos \theta +v_\theta \sin\theta)}.
$$
Therefore, problem \eqref{prob} is converted into the following
problem
\begin{equation}\label{prob-v}
\left\{
\begin{array}{ll}
\displaystyle v_s = 2a\left( \frac{v_\theta \cos \theta - \sin
\theta} {\cos\theta +v_\theta \sin \theta}\right)
\frac{v_{\theta\theta} - v^2_\theta -1}{e^{2v} (\cos\theta +v_\theta
\sin \theta)^2} -1, & |\theta| <\theta_0,\ s> -\frac{1}{2}\log n,\\
v_\theta (-\theta_0,s) = - g_1 (s,v(-\theta_0 ,s)), & s\geqslant -\frac{1}{2}\log n,\\
v_\theta (\theta_0,s) = g_2 (s,v(\theta_0 ,s)), & s\geqslant
-\frac{1}{2}\log n,
\end{array}
\right.
\end{equation}
where
\begin{equation}\label{def:g}
g_i (s,v)=\frac{\sin\theta_0 + k_i \Big(e^{2s}-\frac{1}{n}, e^s e^v
\cos \theta_0\Big)\cos\theta_0}{\cos\theta_0-k_i
\Big(e^{2s}-\frac{1}{n}, e^s e^v \cos \theta_0 \Big) \sin\theta_0 }
\quad (i=1,2).
\end{equation}

\subsection{Gradient bound of $v$}
In a similar way as deriving \eqref{gradient omega} in subsection
\ref{local} one can obtain
\begin{equation}\label{gradient v}
|v_\theta (\theta, s)| \leqslant \Omega_2 (\sigma,\beta) \quad
\mbox{and} \quad \cos \theta +v_\theta \sin \theta \geqslant
\varepsilon_2 (\sigma,\beta)>0
\end{equation}
for $\theta \in [-\theta_0, \theta_0]$ and $s\in [-\frac{1}{2}\log
n, \infty)$.

The second inequality in \eqref{gradient v} implies that, once the
solution $v$ of \eqref{prob-v} is obtained then we can recover it to
the original solution $u$ of \eqref{prob}, since
$$
\frac{\partial }{\partial \theta} X(\theta, v(\theta,s),s) =
\frac{\partial }{\partial \theta} e^s e^{v (\theta,s)} \sin \theta =
e^s e^{v (\theta,s)} (\cos\theta + v_\theta \sin\theta) >0.
$$

\subsection{Bound of $v$}

The local existence result Lemma \ref{lem:local} implies that $v$
exists on $s\in [-\frac{1}{2}\log n, \frac{1}{2}\log (2\tau
+\frac{1}{n})]$. We have changed $u(x,t)$ to a new unknown
$v(\theta,s)$. Similarly, we define $v^\pm $ by $\varphi^\pm$ in the
following way
$$
e^s e^{v^\pm} \cos \theta = \sqrt{2 \Big( e^{2s} -\frac{1}{n}+t^\pm
\Big)}\; \varphi^\pm  \left( \frac{e^s e^{v^\pm} \sin
\theta}{\sqrt{2 \Big( e^{2s} -\frac{1}{n}+t^\pm\Big)}} \right).
$$
By \eqref{bound of u} we have $ e^s e^{v^+} \cos \theta \geqslant
u(e^s e^{v^+} \sin \theta, e^{2s} -\frac{1}{n})$. Noting  $e^s e^{v}
\cos \theta$ $ = u(e^s e^{v} \sin \theta,$ $ e^{2s} -\frac{1}{n})$ we have
$$
(e^{v^+} -e^v)\cos \theta \geqslant u_x \Big(\vartheta, e^{2s}
-\frac{1}{n} \Big) \cdot [(e^{v^+} -e^v)\sin \theta],
$$
where $\vartheta = e^s \sin \theta (\varsigma e^v + (1-\varsigma)
e^{v^+})$ for some $\varsigma \in [0,1]$. Therefore $v^+ \geqslant
v$ by \eqref{regular1}.

On the other hand, by the definition of $v^+$ we have
\begin{equation}\label{bound of v}
e^{v^+} \cos \theta =  \varphi^+ (\cdot) \sqrt{2+ \Big(t^+
-\frac{1}{n}\Big) e^{-2s}} \leqslant \sqrt{2+t^+ /\tau } \; \max
\varphi^{+} \quad \mbox{for } s\geqslant \frac{1}{2}\log \tau.
\end{equation}
So
$$
v(\theta, s) \leqslant v^+ (\theta,s) \leqslant \log \left[
\frac{\sqrt{2+t^+ /\tau } \; \max \varphi^{+}}{\sin \beta}
\right]\quad \mbox{for } \theta\in [-\theta_0, \theta_0],\
s\geqslant \frac{1}{2}\log \tau.
$$

A similar discussion as above shows that
$$
v(\theta, s) \geqslant v^- (\theta,s) \geqslant \log \left[ 2\; \min
\varphi^{-} \right] \quad \mbox{for } \theta\in [-\theta_0,
\theta_0],\ s\geqslant \frac{1}{2}\log \tau.
$$

\begin{remark}
The definition of $v$ depends on $n$, it is not easy to give a
uniform (in $n$) bound for $v$ on $[-\frac{1}{2}\log n, \infty)$, but the
above results show that a uniform bound for $v$  is
possible on $[\frac{1}{2}\log \tau, \infty)$.
\end{remark}

\subsection{Global existence}
Now we consider problem \eqref{prob-v} with initial data
$v(\theta,-\frac{1}{2}\log n)=v_0 (\theta)$, which is defined by
\eqref{v-def-1} at $s=-\frac{1}{2}\log n$. Using the growth bound and gradient bound
in the previous subsections and using the standard theory of parabolic
equations (cf. \cite{Dong, Fri, Lie, LSU}) we can get the following
conclusions.

\begin{lemma}\label{globalv}
    Problem \eqref{prob-v} with initial data $v(\theta,-\frac{1}{2}\log n)
    = v_0 (\theta)$ has a unique, time-global solution $v(\theta,s)
    \in C^{2+\mu, 1 + \mu/2} ([-\theta_0, \theta_0]
    \times [-\frac{1}{2}\log n,\infty))$ and
\begin{equation}\label{est global v}
    \| v(\theta ,s)\|_{C^{2+\mu, 1+ \mu/2}
    ([-\theta_0, \theta_0]\times [\frac{1}{2}\log \tau,\infty)) }
    \leqslant C_0 <\infty,
\end{equation}
where $C_0$ depends on $k_1, k_2, \mu$ and $u_0$ but not on $s$ and $n$.
\end{lemma}

This lemma implies the global existence of $u$.

\begin{lemma}\label{globalu}
    Problem \eqref{prob} and \eqref{ini} has a unique, time-global solution $u(x,t)$.
     Moreover, $u \in C^{2+\mu, 1 + \mu/2} (\overline{Q}_\infty )$,
     where $Q_\infty := \{(x,t)\mid -\xi_1 (t) <x<\xi_2 (t), t>0\}$.
     For any $t_0 > \tau$,
     $$\| u(x,t)\|_{C^{2+\mu, 1+ \mu/2}
     (\overline{Q}_{t_0} \backslash Q_{\tau})} \leqslant C_1 (t_0, k_1, k_2, \mu,u_0,
     \tau ) < \infty.
     $$
\end{lemma}

Indeed, studying the relations between $v$ and $u$ more precisely,
it is not difficult to see that $C_1$ in this lemma can be replaced
by $C_2 \sqrt{t_0} +C_3$ for some $C_2, C_3$ depending on $k_1, k_2, \mu,u_0$
and $\tau$.

\subsection{Existence of self-similar solution}

Since the solution $v$ in Lemma \ref{globalv} is defined for
$s\geqslant -\frac{1}{2}\log n$, we write it as $v_n$. By Cantor's
diagonal argument, one can find a function $V\in C^{2+\mu,
1+\mu/2}([-\theta_0, \theta_0]\times [\frac{1}{2}\log \tau,
\infty))$ and a subsequence $\{n_i\}\subset \{n\}$ such that, as
$i\rightarrow \infty$,
\begin{equation}\label{converge to V}
v_{n_i} (\theta, s) \rightarrow V(\theta,s) \qquad  \mbox{in } \
C^{2, 1}_{\rm loc} \Big([-\theta_0, \theta_0]\times
\Big[\frac{1}{2}\log \tau,\infty \Big) \Big) \ \mbox{topology}.
\end{equation}
Moreover, $V$ satisfies the estimate
\begin{equation}\label{est entire V}
    \| V(\theta ,s)\|_{C^{2+\mu, 1+ \mu/2}
    ([-\theta_0, \theta_0]\times [\frac{1}{2}\log \tau,\infty)) }
    \leqslant C_0 <\infty,
\end{equation}
and $V$ is a solution of
\begin{equation}\label{prob-entire-V}
\left\{
\begin{array}{ll}
\displaystyle v_s = 2a\left( \frac{v_\theta \cos \theta - \sin
\theta} {\cos\theta +v_\theta \sin \theta}\right)
\frac{v_{\theta\theta} - v^2_\theta -1}{e^{2v} (\cos\theta +v_\theta
\sin \theta)^2} -1, & |\theta| <\theta_0,\ s \geqslant \frac{1}{2}\log \tau,\\
v_\theta (-\theta_0,s) = - G_1 (s,v(-\theta_0 ,s)), & s\geqslant \frac{1}{2}\log \tau,\\
v_\theta (\theta_0,s) = G_2 (s,v(\theta_0 ,s)), & s\geqslant
\frac{1}{2}\log \tau,
\end{array}
\right.
\end{equation}
where
\begin{equation}\label{def:H}
G_i (s,v)=\frac{\sin\theta_0 + k_i \Big(e^{2s}, e^s e^v \cos
\theta_0\Big)\cos\theta_0}{\cos\theta_0 -k_i \Big(e^{2s}, e^s e^v
\cos \theta_0 \Big) \sin\theta_0 } \quad (i=1,2).
\end{equation}
Here $G_1$ and $G_2$ are $\log b$-periodic functions in $s$ by \eqref{cond:k}.

Now we use a result in \cite{BPS}.

\begin{lemma}\label{lem:BPS}
Let $u$ be a time-global solution of
\begin{equation}\label{eq-BPS}
\left\{
\begin{array}{ll}
u_t = d(t,x,u,u_x) u_{xx} + f(t,x,u,u_x), & t>0,\ 0<x<1,\\
u_x(i, t) = g^*_i (t, u(i,t)), & i=0,1,\ t>0,
\end{array}
\right.
\end{equation}
where $d,f,g^*_i$ are $C^2$ functions, $T$-periodic in $t$.
If $u(\cdot, t)$ is bounded in $H^2(0,1)$, then there exists a $T$-periodic solution $p$ of \eqref{eq-BPS}
such that $\lim\limits_{t\rightarrow \infty} \|u(\cdot, t)-p(\cdot,
t)\|_{H^2 (0,1)} =0$.
\end{lemma}

By this lemma, problem \eqref{prob-entire-V} has a solution $P(\theta,s)\in H^2(-\theta_0,\theta_0)$,
which is $\log b$-periodic in $s$, and
$$
\lim\limits_{s\to \infty} \|V(\cdot, s) -P(\cdot,s)\|_{H^2(-\theta_0,\theta_0)} =0.
$$
By \eqref{est entire V} we indeed have
\begin{equation}\label{est of P}
\|P (\theta ,s)\|_{C^{2+\mu, 1+ \mu/ 2} ( [-\theta_0, \theta_0]\times\mathbb{R} )}
\leqslant C_0
\end{equation}
and
\begin{equation}\label{V to P}
\lim\limits_{k\to \infty} \|V(\theta,s+k\log b) -P(\theta,s)\|_{C^{2,1}([-\theta_0,\theta_0]\times
[0,\log b])} =0.
\end{equation}

We now recover $P$ to a solution of $\eqref{prob}$, the corresponding change
of variables should be the limiting version as $n\rightarrow \infty$
of \eqref{new-variable-v} and \eqref{inverse-map-v}, or for $s\in
(-\infty, \infty)$. Using these variables we define $U$ by
\begin{equation}\label{P-def-U}
e^s e^P \cos\theta = U(e^s e^P \sin \theta, e^{2s})\quad \mbox{for }
\theta\in [-\theta_0, \theta_0],\ s\in (-\infty, \infty).
\end{equation}
By \eqref{gradient v}, \eqref{converge to V} and \eqref{V to P} we have
$$
\frac{\partial}{\partial \theta} X ( \theta, P(\theta,s), s)=
\frac{\partial}{\partial \theta} \left( e^s e^{P(\theta,s)}
\sin\theta\right) = e^s e^{P(\theta,s)} (P_\theta \sin\theta +
\cos\theta) >0,
$$
so the function $U(x,t)$ is well-defined for all $t>0$. Moreover, by
the definition of $U$ and the periodicity of $P$ we have
\begin{eqnarray*}
b U(e^s e^{P(\theta,s)} \sin \theta, e^{2s}) & = &  b e^s
e^{P(\theta, s)} \cos\theta  =  e^{s +\log b} e^{P(\theta, s+\log
b)} \cos\theta \\
& = & U(e^{s+\log b} e^{P(\theta, s+\log b)} \sin \theta, e^{2s+2\log b} ) \\
& = & U(b e^s e^{P(\theta, s)} \sin \theta, b^2 e^{2s} ).
\end{eqnarray*}
Hence
$$
bU(x,t) = U(bx, b^2 t )\quad \mbox{for } -\Xi_1 (t)\leqslant x
\leqslant \Xi_2 (t),\ t>0,
$$
where $-\Xi_1(t)$ and $\Xi_2(t)$ are the $x$-coordinate of the end
points of the graph of $U$. Since $t=e^{2s} $ we have
$$
\Xi_2 (t) = e^s e^{P(\theta_0, s)} \sin \theta_0 = \sqrt{t} e^{P (
\theta_0, \frac{1}{2}\log t )} \sin\theta_0\quad \mbox{for } t>0.
$$
So
$$
\Xi_2 (b^2 t) = b \sqrt{t} e^{P (\theta_0, \frac{1}{2}\log t +\log
b)}  \sin\theta_0 = b \sqrt{t} e^{P(\theta_0, \frac{1}{2}\log t )}
 \sin\theta_0 = b\Xi_2 (t)\quad \mbox{for } t>0.
$$
Similarly we have $\Xi_1 (b^2 t) = b\Xi_1 (t)$ for $t>0$.
Consequently, we obtain a discrete expanding self-similar solution of \eqref{prob} and
this proves the existence part of Theorem \ref{thm:1}.

\subsection{Asymptotic stability}
In this subsection we assume that $k_i (t,u) \equiv k_i (u)\ (i=1,2)$ and
prove the asymptotic stability result in Theorem \ref{thm:1}:
$$
d_{\mathcal{H}} (\Gamma(t), \gamma(t)) \leqslant C t^{-1/2}, \quad t\to \infty,
$$
where $\Gamma(t)$ is the graph of the discrete expanding self-similar solution $U$,
$\gamma(t)$ is the graph of any solution $u$ of \eqref{prob} with some initial data,
and $d_\mathcal{H}$ is the Hausdorff distance.

When $k_1$ and $k_2$ are constants satisfying \eqref{cond:expand},
in \cite{CGK, CGL, Koh}, the authors proved similar results for problem
\eqref{prob0} by constructing precise lower and upper solutions.
We will use the change of variables and the a priori estimates but do not
construct lower and upper solutions.
So our approach is different from those in \cite{CGK, CGL, Koh}.

For any initial data $u_0\in C^{2+\mu}_{\rm ad}$, there exists
$t_1 >0$ such that
$$
0= U(\cdot,0)\leqslant u_0 (\cdot) \leqslant U(\cdot,t_1).
$$
By comparison principle Lemma \ref{comparison} we have
\begin{equation}\label{U-ut-Ut1}
U(\cdot,t)\leqslant u (\cdot,t) \leqslant U(\cdot,t+t_1),\quad t>0.
\end{equation}

For any given $t>0$ and any $\bar{x}$ with $(\bar{x}, u(\bar{x},t)) \in \overline{S}$,
denote $\bar{\theta} = \arctan \frac{\bar{x}}{u(\bar{x},t)}$. The line $\theta =
\bar{\theta}$ contacts $\Gamma(t)$ (resp. $\Gamma(t+t_1)$) at exactly one point
$A$ (resp. $A_1$). Denote the $x$-coordinate of $A$ (resp. $A_1$) by
$$
Z = Z(\bar{\theta}) = Z (\bar{x},t) \qquad (\mbox{resp. } \ Z_1 = Z_1
(\bar{\theta}) = Z_1 (\bar{x},t)\; ).
$$
It follows from \eqref{U-ut-Ut1} and the equivalence in \eqref{op-1} that
\begin{equation}\label{Z1<u<Z2}
[Z^2 +U^2 (Z, t)]^{1/2} \leqslant [\bar{x}^2 + u^2 (\bar{x}, t)]^{1/2}
\leqslant [Z_1^2 +U^2 (Z_1, t+t_1)]^{1/2} .
\end{equation}
Set $s =\frac{1}{2} \log t$ and $s(t_1)  = \frac{1}{2} \log (t+t_1)$.
By the definition of $U$ in \eqref{P-def-U} we have
$$
Z = e^s e^{P(\bar{\theta}, s)} \sin\bar{\theta},\quad \ \
U(Z, t) = e^s e^{P(\bar{\theta}, s)} \cos\bar{\theta},
$$
and
$$
Z_1 = e^{s(t_1)} e^{P(\bar{\theta}, s(t_1))} \sin\bar{\theta},\quad \ \
U(Z_1, t+t_1) = e^{s(t_1)} e^{P(\bar{\theta}, s(t_1))} \cos\bar{\theta},
$$
where $P$ is that in \eqref{est of P}. Since
$$
s(t_1) - s = \frac{t_1}{2t } + o\Big(\frac{1}{t}\Big)\ \ (\mbox{as }
t \to \infty),
$$
by \eqref{est of P} we have
\begin{eqnarray*}
[Z_1^2 +U^2 (Z_1, t+t_1)]^{1/2} - [Z^2 +U^2 (Z, t)]^{1/2}
 & = & e^{s(t_1)} e^{P(\bar{\theta}, s(t_1))} - e^{s} e^{P(\bar{\theta}, s)}\\
&\leqslant & C_1 t^{-1/2}, \quad  \mbox{as } t\to \infty,
\end{eqnarray*}
where $C_1$ depends on $t_1$ and $C_0$ in \eqref{est of P}, but not on $\bar{\theta}$
and $t$.
Therefore, we have
$$
d_\mathcal{H} (\Gamma(t), \gamma(t)) \leqslant
d_\mathcal{H} (\Gamma(t), \Gamma(t+t_1)) \leqslant
C_1 t^{-1/2}\quad \mbox{ as } t\to \infty.
$$

\subsection{Uniqueness of self-similar solutions}
In this subsection we still assume $k_i (t,u) \equiv k_i (u)\ (i=1,2)$ and
to prove the uniqueness conclusion in Theorem \ref{thm:1}. The uniqueness
for general $k_i (t,u)$ is still open.

We begin with choosing a convex initial data $u_0$, that is,
$$
a(u_{0x}) u_{0xx} \geqslant \epsilon
$$
for some $\epsilon>0$. Such a choice is possible. For example, draw a line $\ell_1$ from $A_1
:=(-1,\tan\beta)$ with slope $-k_1 (\tan\beta)$. Denote the contacting
point between $\ell_1$ and $\partial_2 S$ by $A'_2$, then
$A'_2:= (x',x'\tan\beta)$ with
$$
x' := \frac{\tan\beta -k_1 (\tan\beta)}{\tan\beta +k_1 (\tan\beta)}.
$$
Choose $A_2 := (x'+\epsilon', (x'+\epsilon')\tan\beta)\in \partial_2 S$ for
some small $\epsilon'>0$, then $A_2$ is above $A'_2$. Draw a line $\ell_2$
from $A_2$ with slope $k_2 ((x'+\epsilon')\tan\beta)$. Since
$$
k_2 ((x'+\epsilon')\tan\beta) \geqslant k^0_2 > -k^0_1 \geqslant
-k_1 (\tan\beta),
$$
$\ell_2$ must contact $\ell_1$ at some point $A_3\in S$ provided
$\epsilon'>0$ is small enough.
Now we smoothen $\overline{A_1 A_3}+\overline{A_3 A_2}$ such that
the smoothened curve $\mathcal{C}$ is strictly convex, it is tangent to $\overline{A_1 A_3}
$ at $A_1$, tangent to $\overline{A_3 A_2}$ at $A_2$. Now the corresponding
function $u_0$ of $\mathcal{C}$ is a desired initial data.

Let $u$ be the solution of \eqref{prob} with the above constructed initial
data $u_0$. Denote $\eta = u_t$. Differentiating the problem \eqref{prob} by $t$ we have
$$
\left\{
\begin{array}{l}
 \eta_t = a(u_x) \eta_{xx} + a'(u_x) u_{xx} \eta_{x},\quad -\xi_1 (t)<x<\xi_2(t),\ t>0,\\
 \eta_x (-\xi_1 (t),t) = f_1 (t)\eta (-\xi_1 (t),t),\ \
 \eta_x (\xi_2 (t),t) = f_2 (t)\eta (\xi_2 (t),t), \ \ \quad \ t>0,\\
 \eta(x,0) = a(u_{0x}) u_{0xx} \geqslant \epsilon >0,
\end{array}
\right.
$$
where $f_1$ and $f_2$ are continuous functions. Maximum principle implies
that $\eta = u_t >0$ for $t>0$. In a similar way as in the previous subsection
3.6 we can get a discrete expanding self-similar solution $(U, \Xi_1, \Xi_2)$.
Moreover, using the same notions as above we have $1+v_s >0 $ and so
$1+V_s \geqslant 0$ for $s\geqslant \frac{1}{2}\log \tau$. Using \eqref{V to P}
one has $1+P_s \geqslant 0$ for all $s\in \R$, this implies that $U_t \geqslant 0$. Finally, the
strong maximum principle implies that $U_t >0$ for all $t>0$, and so
$\Xi_{1t}(t)>0,\ \Xi_{2t}(t)>0$ for $t>0$.

Suppose that $(U^*, \Xi^*_1, \Xi^*_2)$ is another discrete self-similar solution of \eqref{prob}.
We want to prove that $U(x,t)\equiv U^*(x,t)$. Otherwise, either
\begin{equation}\label{contra 1}
U^*(\cdot,t) \leqslant U(\cdot,t) \ \ \ \mbox{and}\ \ \  U(x,t)\not\equiv U^*(x,t)
\ \quad \mbox{for all } t>0,
\end{equation}
or
\begin{equation}\label{contra 2}
U(\cdot,t) \leqslant U^*(\cdot,t) \ \ \ \mbox{and}\ \ \  U(x,t)\not\equiv U^*(x,t)
\ \quad \mbox{for all } t>0,
\end{equation}
holds.
We now derive a contradiction from \eqref{contra 1}.

Since $U(x,t) \to 0$ as $t \to 0$, for any $t_0 >0$, there exists
$\tau\in (0,t_0)$ such that $ U(\cdot ,t_0-\tau) \ll U^* (\cdot,t_0)$.
Set
$$
\tau(t_0) := \inf\{ \tau \mid U(\cdot,t_0 -\tau) \ll U^* (\cdot,t_0)\}.
$$
Then either
\begin{itemize}
\item[\rm (a)] $U(x,t_0 -\tau(t_0)) \equiv U^* (x,t_0)$, or

\item[\rm (b)] $U(x,t_0 -\tau(t_0)) \not\equiv U^* (x,t_0)$ and $U(\cdot ,t_0 -\tau(t_0)) \preceq U^* (\cdot,t_0)$
\end{itemize}
holds. By \eqref{contra 1} and by $U_t >0$ it is easily seen that $\tau(t_0) >0$, and so $\tau(t_0) \in (0,t_0)$.

We first show that (b) is impossible. Otherwise,  by Lemma \ref{comparison} we have
$$
U(\cdot, t_0 -\tau(t_0)+ (b^2 -1) t_0 ) \ll  U^* (\cdot, t_0 +(b^2 -1) t_0) = U^*(\cdot, b^2 t_0).
$$
On the other hand, by the self-similarity we have
$$
U(bx, b^2 (t_0 - \tau(t_0))) \equiv bU(x,t_0 -\tau(t_0))\preceq b U^*(x,t_0)\equiv
U^*(bx,b^2 t_0).
$$
Combining these inequalities with the fact that $U_t >0$ we have
$$
b^2 t_0 -\tau(t_0) < b^2 (t_0 - \tau(t_0)).
$$
This is a contradiction.

Now, (a) implies that,  for any $t>0$, there exists $\tau(t) \in (0,t)$ such that
$U(x,t -\tau(t)) \equiv U^* (x,t)$. Hence, for any $t,s>0$ we have
$$
U(x,t-\tau(t) +s)\equiv U^*(x,t +s) \equiv U(x,t+s -\tau(t+s)).
$$
So $\tau(t+s)=\tau(t)$ for all $t,s>0$. Thus, $\tau(t)$ is a constant $\tau^* > 0$
and $U(x,t-\tau^*) \equiv U^*(x,t)$ for all $t>0$. Taking limits as $t\to \tau^*+0$ we have
$U^*(x,\tau^*) = U(x, 0) = 0$, and so $\tau^* =0$, a contradiction.

The above discussion shows that \eqref{contra 1} is impossible.
In a similar way, one can show that \eqref{contra 2} is also impossible.
So $U(x,t)\equiv U^*(x,t)$, and $(U, \Xi_1, \Xi_2)$ is the unique
discrete expanding self-similar solution
of \eqref{prob}. This proves the uniqueness result in Theorem \ref{thm:1}.


\section{Shrinking self-similar solutions}
\par
In this section we always assume that \eqref{cond:shrink} holds.

\subsection{Classical shrinking/backward self-similar solutions}
First, recall the classical self-similar solutions of \eqref{prob0}.
For any $\gamma_1, \gamma_2 \in \R$, consider the problem
\begin{equation}\label{eq-psi}
\left\{
\begin{array}{l}
a(\psi'(z))\psi''(z) = z\psi'(z) -\psi (z),\ \ \ z\in \R,\\
\psi'(-q_1) = -\gamma_1, \ \psi(-q_1) = q_1 \tan\beta,\\
\psi'(q_2) = \gamma_2, \ \psi(q_2) = q_2\tan\beta.
\end{array}
\right.
\end{equation}
In \cite{GH, GK}, the authors obtained the following result.

\begin{lemma}
For any $\gamma_1, \gamma_2$ satisfying $\gamma_1 + \gamma_2 <0$,
there exists a pair $q_1, q_2 >0$ such that problem \eqref{eq-psi}
has solution $\psi(z;\gamma_1, \gamma_2)$, which is positive on
$[-q_1, q_2]$.
\end{lemma}

It is easily seen that the function
$$
\sqrt{2(T-t)}\; \psi\left( \frac{x}{ \sqrt{2(T-t)} }\right) \quad
\mbox{for } - \zeta_1(t) < x < \zeta_2(t),\ t>0,
$$
with $\zeta_i (t) = q_i \sqrt{2(T-t)}\ (i=1,2)$ is a classical
shrinking/backward self-similar solution of \eqref{prob0}.

We use $(\psi^- (z), q^-_1, q^-_2)$ to denote the solution of
\eqref{eq-psi} with $\gamma_i = k^0_i\ (i=1,2)$, use $(\psi^+ (z),
q^+_1, q^+_2)$ to denote the solution of \eqref{eq-psi} with
$\gamma_i = K^0_i\ (i=1,2)$, where $k^0_i$ and $K^0_i$ are those in \eqref{cond:k2}.

By \eqref{cond:shrink}, $K^0_1 +K^0_2 <0$ and so $(\psi^\pm )'' <0$.
For any $T_0 >0$, the function $\sqrt{2(T_0 -t)}\; \psi^-
(x/\sqrt{2(T_0 -t)})$ and the function $\sqrt{2(T_0 -t)}\; \psi^+
(x/\sqrt{2(T_0 -t)})$ (both are shrinking/backward self-similar
solutions of \eqref{prob0}) are lower and upper solutions of
\eqref{prob}, respectively. Since the initial data $u_0 >0$, there
exists $T^+ , T^- >0$ such that
\begin{equation}\label{less and equal}
\sqrt{2T^- }\; \psi^-  \left({\frac {\cdot} {\sqrt{2T^-}}}\right)
\preceq u_0 (\cdot) \preceq \sqrt{2T^+ }\; \psi^+ \left({\frac {\cdot}
{\sqrt{2T^+ }}}\right).
\end{equation}
Comparison principle implies that
\begin{equation}\label{bound of u shrink}
\sqrt{2(T^- -t)}\; \psi^- \left({\frac {\cdot} {\sqrt{2(T^- -
t)}}}\right) \ll u(\cdot,t) \ll  \sqrt{2(T^+ -t)}\; \psi^+
\left({\frac {\cdot} {\sqrt{2(T^+ - t)}}}\right)
\end{equation}
on the time interval where these three functions are defined.

\begin{lemma}\label{lem:T^-/T^+}
Let $\psi^\pm$ and $T^\pm$ be as in \eqref{less and equal}. Then
$T^+ > T^- \geqslant \delta T^+$, where
$$
\delta = \delta(\beta, \sigma, k^0_i, K^0_i) :=
\frac{\sigma^4}{(2\tan\beta -\sigma)^4 } \cdot
\frac{(q^+_2)^2}{(q^-_2)^2} >0.
$$
\end{lemma}

\begin{proof}
We first prove $T^+ > T^-$. The areas $D^\pm (t)$ of
the regions enclosed by the graph of $\sqrt{2(T^\pm -t)}\; \psi^\pm
(x/ \sqrt{2(T^\pm -t)} )$, $\partial_1 S$ and $\partial_2 S$ are
given by
$$
D^\pm (t) = \int_{-\zeta_1 (t)}^{\zeta_2(t)} \sqrt{2(T^\pm -t)}\;
\psi^\pm \left({\frac {x} {\sqrt{2(T^\pm -t)}}}\right) dx -
\frac{1}{2}[\zeta^2_1  (t) + \zeta_2^2 (t)]\tan\beta.
$$
A simple computation shows that
$$
(D^+)'(t) = \int_{-K^0_1}^{K^0_2} a(p) dp,\quad (D^-)'(t) =
\int_{-k^0_1}^{k^0_2} a(p) dp.
$$
Since $D^+(T^+) =D^- (T^-)=0$ we have
\begin{equation}\label{area est}
D^+ (0) = \int_0^{T^+} dt\int_{K^0_2}^{-K^0_1} a(p) dp,\quad D^- (0)
= \int_0^{T^-} dt\int_{k^0_2}^{-k^0_1} a(p) dp.
\end{equation}
By \eqref{cond:shrink} we have $k^0_2 <K^0_2 < -K^0_1 < -k^0_1$,
and by \eqref{less and equal} we have $D^+ (0) \geqslant D^-(0)$. Therefore,
we have $T^+ >T^-$ by \eqref{area est}.

Next we prove $T^- \geqslant \delta T^+$. For $i=1,2$, denote
$Q^+_i$ (resp. $Q^-_i$) the end points of the graph of
$\sqrt{2T^+}\; \psi^+ (x/ \sqrt{2T^+} )$ (resp. the graph of
$\sqrt{2T^-}\; \psi^- (x/ {\sqrt{2T^-}} ))$ on $\partial_i S$,
respectively.

Connecting $Q^+_1$ and $Q^+_2$ we get a line segment
$\overline{Q^+_1 Q^+_2}$. It is below the graph of $\sqrt{2T^+}\;
\psi^+ ( x/ \sqrt{2T^+})$ since $(\psi^+)'' <0$. Draw a line from
$Q^+_1$ (resp. $Q^+_2$) with slope $-\tan\beta +\sigma$ (resp.
$\tan\beta -\sigma$). Assume that it contacts $\partial_2 S$ (resp.
$\partial_1 S$) at $A_2$ (resp. $A_1$).

By \eqref{less and equal} the graph of $u_0$ contacts the line
segment $\overline{Q^+_1 Q^+_2}$. Since $|u_{0x}|\leqslant \tan\beta
-\sigma$, we see that the graph of $u_0$ is above the line segment
$\overline{A_1 A_2}$.

If for some $T_0 >0$, the graph of $\sqrt{2T_0}\; \psi^- (x/
\sqrt{2T_0})$ is tangent to $\overline{A_1 A_2}$ from above, then (note
$|(\psi^-)'|\leqslant \tan\beta -\sigma$) the graph of $\sqrt{2T_0}\;
\psi^- (x/ \sqrt{2T_0} )$ lies above the line segment $\overline{B_1
B_2}$, where $B_1$ (resp. $B_2$) is the contacting point between the
line passing $A_2$ (resp. $A_1$) with slope $\tan\beta -\sigma$
(resp. $-\tan\beta +\sigma$) and the left boundary $\partial_1 S$
(resp. right boundary $\partial_2 S$). Therefore, $Q^-_i$ are above
$B_i$ for $i=1,2$.

Using the coordinates of $Q^+_i = ((-1)^i r_i\cos\beta, r_i
\sin\beta )$, where
$$
r_i = \frac{\sqrt{2T^+} q^+_i }{\cos\beta} \quad (i=1,2),
$$
one can easily calculate the coordinates of $B_1$ and $B_2$:
$$
B_i =\Big( \frac{(-1)^i \sigma^2 r_i \cos\beta}{(2\tan\beta
-\sigma)^2 }, \ \frac{\sigma^2 r_i \sin\beta}{(2\tan\beta -\sigma)^2
}\Big)\quad (i=1,2).
$$
The fact that $Q^-_2$ is above $B_2$ implies that
$$
\sqrt{2T^-} q^-_2 \geqslant \frac{\sigma^2 \cos\beta}{(2\tan\beta
-\sigma)^2 } \cdot \frac{\sqrt{2T^+} q^+_2 }{\cos\beta}.
$$
So
$$
\frac{T^-}{T^+} \geqslant \left[ \frac{\sigma^2
q^+_2}{q^-_2(2\tan\beta -\sigma)^2} \right]^2.
$$
This proves the lemma.
\end{proof}

\subsection{Shrinking time for solutions of \eqref{prob} and \eqref{ini}}
In this subsection we consider the shrinking time for the
solution $u=u(x,t;u_0)$ of \eqref{prob} with initial data $u_0$.
We give two results. The first one is about the shrinking
time of $u=u(x,t;u_0)$ for any given $u_0$, the second one is about
the existence of $u_0$ for given shrinking time $T$.

\begin{lemma}\label{lem:exist converge}
Let $u_0\in C^{2+\mu}_{\rm ad}\ (\mu \in (0,1))$ be an initial data. If it
satisfies \eqref{less and equal}, then there
exists $T\in (T^-,T^+)$ such that the solution $u(x,t;u_0)$ of
\eqref{prob} and \eqref{ini} exists on time interval $[0,T)$, and
\begin{equation}\label{convergence}
\|u(\cdot ,t)\|_{L^\infty} \rightarrow 0, \quad \xi_1 (t)\rightarrow
0, \quad \xi_2 (t)\rightarrow 0\quad \mbox{as } t\rightarrow T.
\end{equation}
\end{lemma}

\begin{proof}
We use polar coordinates $x=r\sin\theta,\
y=r\cos\theta$ for $\theta\in [-\theta_0, \theta_0]$, that is,
$$
r\cos\theta = u(r\sin \theta, t)
$$
defines an implicit function $r=r(\theta, t)$ by \eqref{regular1}.
Problem \eqref{prob} is then converted into
\begin{equation}\label{prob-r}
\left\{
\begin{array}{ll}
\displaystyle r_t = a\left( \frac{r_\theta \cos\theta - r
\sin\theta} {r_\theta \sin\theta +r\cos\theta }
\right)\frac{rr_{\theta\theta} -2r^2_\theta -r^2 }{r(r_\theta
\sin\theta +r\cos\theta)^2}, &
\theta\in [-\theta_0, \theta_0],\ t>0,\\
r_\theta (-\theta_0, t) = - \tilde{h}_1 (t,r(-\theta_0, t)),&
t\geqslant 0,\\
r_\theta (\theta_0, t) = \tilde{h}_2 (t,r(\theta_0, t)),&
t\geqslant 0,
\end{array}
\right.
\end{equation}
where
$$
\tilde{h}_i (t,r) := \frac{\sin\theta_0 + k_i(t,r\cos\theta_0) \cos\theta_0
} { \cos\theta_0 - k_i(t,r\cos\theta_0) \sin\theta_0}r\quad (i=1,2).
$$

By \eqref{bound of u shrink} we have $r\cos\theta =
u(r\sin\theta, t)\leqslant \sqrt{2T^+} \max\psi^+$. So $
0\leqslant r(\theta,t)\leqslant \sqrt{2T^+} \max\psi^+ /\sin\beta$.
By \eqref{gradient bound} and \eqref{regular1} we have
$$
|r_\theta| =\left| \frac{\sin\theta +u_x \cos\theta}{\cos\theta -u_x
\sin\theta} r \right|\leqslant \frac{1 +\tan\beta -\sigma}{\sigma
\cos\beta } \cdot \frac{\sqrt{2T^+} \max\psi^+}{\sin\beta}.
$$
Thus the standard a priori estimates (cf. \cite{Dong, Fri, Lie, LSU})
show that the solution of \eqref{prob-r} with initial data $r(\theta,0)$
will not develop singularity till $\min r(\cdot, t) \rightarrow 0$
as $t\rightarrow T$ for some $T>0$. Moreover, $T^- < T < T^+ $ follows from
\eqref{bound of u shrink} and Lemma \ref{lem:T^-/T^+}.

Now we prove \eqref{convergence}. If $u(0,T+0)=0$ but $u(\bar{x}, T+0)>0$
for some $\bar{x}\not= 0$, then there exists $\hat{x}$ lying between $0$ and $\bar{x}$
such that
$$
|u_x (\hat{x}, T+0)| =
\frac{|u(\bar{x},T+0)-u(0,T+0)|} {|\bar{x}|}\geqslant \tan\beta.
$$
This contradicts Lemma \ref{lem:gradient bound} and so the first limit in
\eqref{convergence} holds. The last two limits in \eqref{convergence}
follow from the first one. This proves the lemma.
\end{proof}

\begin{lemma}\label{exist ini}
For any given $T>0$, there exists an initial data
$u_0$ such that the solution $u(x,t;u_0)$ of \eqref{prob} and \eqref{ini}
exists on $[0,T)$ and it shrinks to $0$ just at time $T$.
\end{lemma}

\begin{proof}
Choose two initial data $u_{01}, u_{02} \in C^{2+\mu}_{\rm ad}$.
Moreover, $u_{01}$ is chosen so large such that \eqref{less and equal}
holds for $u_0 = u_{01}$  and for some $T^- >T$. By Lemma \ref{lem:exist converge},
the solution $u(x,t;u_{01})$ shrinks to $0$ as $t\rightarrow \widetilde{T}_1$
for some $\widetilde{T}_1 >T^- >T$.
On the other hand, we choose $u_{02}$ small such that \eqref{less and equal}
holds for $u_0 = u_{02}$ and for some $T^+ < T$. By Lemma \ref{lem:exist converge} again, the
solution $u(x,t;u_{02})$ shrinks to $0$ as $t\rightarrow \widetilde{T}_2$ for some
$\widetilde{T}_2 <T^+ <T$.

Now we modify the initial data from $u_{02}$ to $u_{01}$ little by little
such that the modified initial data is still in $C^{2+\mu}_{\rm ad}$.
Since the solution $u(x,t;u_0)$ of \eqref{prob} and \eqref{ini} depends on
the initial data $u_0$ continuously, we finally have an initial data
$u_0$ such that $u(x,t;u_0)$ shrinks to $0$ at time $T\in (\widetilde{T}_2,
\widetilde{T}_1)$.
\end{proof}

In the following, we fix $T>0$ and choose the initial data
as in Lemma \ref{exist ini}.

\subsection{Change of variables}
Lemma \ref{lem:exist converge} gives the existence and boundedness of
$r$ (and so, of $u$), but the time interval is finite: $[0,T)$. So Lemma
\ref{lem:BPS} can not be applied to give a periodic solution.
To get a shrinking self-similar solution of \eqref{prob}, we
introduce new coordinates.  Set
\begin{equation}\label{new-variable-shrink}
\left\{
\begin{array}{ll}
\displaystyle \theta = \arctan \frac{x}{y}, & (x,y)\in
\overline{S}\backslash
\{0\}, \\
\displaystyle \rho= -\frac{1}{2} \log \frac{x^2 +y^2}{T-t} , &
(x,y)\in
\overline{S} \backslash \{0\},\ 0\leqslant t <T,\\
\displaystyle s= -\frac{1}{2} \log (T-t),& 0\leqslant t <T.
\end{array}
\right.
\end{equation}
The inverse map is
\begin{equation}\label{inverse-map-2}
\left\{
\begin{array}{ll}
x= e^{-s} e^{-\rho} \sin \theta,& (\theta, \rho)\in \overline{D},\ s\in [-\frac{1}{2}\log T, \infty),\\
y=e^{-s} e^{-\rho} \cos \theta,& (\theta, \rho)\in \overline{D},\ s\in [-\frac{1}{2}\log T, \infty),\\
t = T- e^{-2s}, & s\in [-\frac{1}{2}\log T, \infty).
\end{array}
\right.
\end{equation}
A similar discussion as in subsections 2.4 and 3.2 shows that
in these new variables, the original function
$y=u(x,t)$ is converted into a new function $\rho = w(\theta,s)$.
Differentiating the expression
\begin{equation}\label{def of w}
e^{-s} e^{-w(\theta,s)} \cos \theta = u(e^{-s} e^{-w(\theta,s)} \sin
\theta, T-e^{-2s})
\end{equation}
twice by $\theta$ and once by $s$ we obtain
\begin{equation}\label{exp u w}
 u_x = \frac{w_\theta
\cos \theta + \sin\theta}{w_\theta \sin \theta - \cos \theta}, \quad
u_{xx} =  \frac{e^s e^w(w_{\theta\theta} +w^2_\theta +1)}{(w_\theta
\sin \theta - \cos \theta)^3}, \quad u_t = \frac{e^s(1+w_s)}{2e^w
(w_\theta \sin \theta - \cos \theta)} .
\end{equation}
Therefore, problem \eqref{prob} is converted into the following
problem
\begin{equation}\label{prob-w}
\left\{
\begin{array}{l}
\displaystyle w_s = 2e^{2w} a\left( \frac{w_\theta \cos \theta
+\sin\theta} {w_\theta \sin \theta - \cos \theta}\right) \frac{
w_{\theta\theta} + w^2_\theta +1} {(w_\theta \sin \theta -
\cos\theta )^2 } -1, \\
 \ \ \ \ \ \ \ \ \ \ -\theta_0 <\theta <\theta_0,\
s\in [-\frac{1}{2}\log T, \infty),\\
w_\theta (-\theta_0,s) = \tilde{g}_1 (s,w(-\theta_0 ,s)), \quad
s\in [-\frac{1}{2}\log T, \infty),\\
w_\theta (\theta_0,s) = -\tilde{g}_2 (s,w(\theta_0 ,s)), \quad s\in
[-\frac{1}{2}\log T, \infty),
\end{array}
\right.
\end{equation}
where
\begin{equation}\label{def:tilde g}
\tilde{g}_i (s,w)=\frac{\sin \theta_0 +k_i (T-e^{-2s}, e^{-s} e^{-w} \cos \theta_0)
\cos \theta_0 }{\cos \theta_0 - k_i (T-e^{-2s},
e^{-s}e^{-w} \cos \theta_0) \sin \theta_0} \quad (i=1,2).
\end{equation}

\subsection{Bound of $w$}
We derive the boundedness of $w$ in a series time intervals:
$[0,\delta^2 T],\ [\delta^2 T, T-(1-\delta^2)^2 T],\
[T-(1-\delta^2)^2 T, T-(1-\delta^2)^3 T],\cdots$, where $\delta\in (0,1)$ is
that in Lemma \ref{lem:T^-/T^+}.

In the first step, we choose $\psi^\pm$ and $T^\pm$ as in
\eqref{less and equal}, and consider \eqref{prob} on time interval
$t\in [0,\delta^2 T]$ (note that $\delta^2 T < \delta^2 T^+
\leqslant \delta T^- < \delta T < \delta T^+ \leqslant T^-$), or,
equivalently, consider \eqref{prob-w} on time-interval $s\in
[-\frac{1}{2}\log T, -\frac{1}{2}\log T
-\frac{1}{2}\log(1-\delta^2)]$. In this period,
$$
e^{-2s} = T- t \geqslant (1-\delta^2)T  \quad \Rightarrow \quad
e^{2s} T \leqslant \frac{1}{1-\delta^2}.
$$
Thus,
\begin{eqnarray*}
e^{2s}(T^+ -t) & \leqslant & e^{2s} [(T^+ -T^-) + (T-t)]  \leqslant
e^{2s}
\left[ \frac{1-\delta}{\delta} T^- + (T-t) \right] \\
& \leqslant & \frac{1-\delta}{\delta} e^{2s} T + 1
\leqslant \delta_1 := 1+ \frac{1}{\delta(1+\delta)}.
\end{eqnarray*}

By \eqref{bound of u shrink} we have, for $t\in [0,\delta^2 T]$ or $
s\in [-\frac{1}{2}\log T, -\frac{1}{2}\log T -\frac{1}{2} \log
(1-\delta^2)]$,
$$
e^{-s}e^{-w} \cos \theta = u(\cdot, T-e^{-2s}) \leqslant \sqrt{2(T^+
-t)} \max \psi^+.
$$
So
$$
 e^{-w} \leqslant \frac{\max \psi^+}{\cos\theta_0} e^{s}\sqrt{2(T^+
-t)} \leqslant \frac{\max \psi^+}{\cos\theta_0} \sqrt{2\delta_1}.
$$

On the other hand, in the same time interval $t\in [0,\delta^2 T]$
we have
$$
\frac{T-T^-}{T-t} \leqslant \frac{T-T^-}{T- \delta^2 T}  \leqslant
\frac{T- \delta T}{T- \delta^2 T} = \frac{1 }{1 + \delta }.
$$
So
$$
e^{2s} (T^- -t) = \frac{(T^- - T)+(T-t)}{T-t} \geqslant 1-
\frac{1}{1+\delta} = \frac{\delta}{1+\delta}.
$$
Thus by
$$
e^{-s}e^{-w} \cos \theta = u(\cdot, T-e^{-2s}) \geqslant \sqrt{2(T^-
-t)} \min \psi^-,
$$
we have
$$
 e^{-w} \geqslant \min \psi^-  e^{s}\sqrt{2(T^-
-t)} \geqslant \min \psi^- \sqrt{\frac{2\delta}{1+\delta}}.
$$
Therefore we obtain the bound of $w$ for $s\in [-\frac{1}{2}\log T,
-\frac{1}{2}\log T -\frac{1}{2}\log(1-\delta^2)]$:
\begin{equation}\label{bound of w}
-\log \left[ \frac{\max \psi^+}{\cos\theta_0} \sqrt{2\delta_1}
\right] \leqslant w \leqslant -\log \left[ \min \psi^-
\sqrt{\frac{2\delta}{1+\delta}}\right].
\end{equation}
Note that the lower and upper bounds do not depend on $s$.

Take a another pair $T^+_*, T^-_* >0$ such that
$$
\sqrt{2T^-_* }\; \psi^- \left(\frac{\cdot}{\sqrt{2T^-_* }}
\right) \preceq u(\cdot,\delta^2 T) \preceq \sqrt{2T^+_*}\; \psi^+
\left( \frac{\cdot}{\sqrt{2T^+_* }} \right).
$$
Lemma \ref{lem:T^-/T^+} implies that $T^+_* >
T^-_* \geqslant \delta T^+_*$. From time $\delta^2 T$, $u$
will shrink to $0$ in time $T_2 := (1-\delta^2)T$, we consider
another time interval: $t\in [\delta^2 T, \delta^2 T +\delta^2 T_2] =
[\delta^2 T, T-(1-\delta^2)^2 T]$, or $s\in [-\frac{1}{2}\log T -\frac{1}{2}
\log (1-\delta^2),  -\frac{1}{2}\log T - \log (1-\delta^2)]$. Replacing
$T^\pm$ by $T^\pm_*$ in the above discussion we see that
\eqref{bound of w} holds on this time interval.

Repeat such processes infinite times we obtain the estimate
\eqref{bound of w} for $w$ on $[0,T)$.

\subsection{A priori estimate for $w$}
The gradient bound of $w$ is similar as that for $\omega$ in subsection 2.4
and that for $v$ in subsection 3.3. Using the standard theory
of parabolic equations (cf. \cite{Dong, Fri, Lie, LSU}) we can get the
following conclusions.

\begin{lemma}\label{globalw}
    Problem \eqref{prob-w} with initial data $w(\theta, -\frac{1}{2}\log T)$
    (which is defined by \eqref{def of w} at $s=-\frac{1}{2}\log T$) has a unique, time-global
    solution $w(\theta,s)     \in C^{2+\mu, 1 + \mu/2} $ $ ([-\theta_0, \theta_0]
    \times [-\frac{1}{2}\log T, \infty))$ and
    $$
    \| w(\theta ,s)\|_{C^{2+\mu, 1+ \mu/ 2}
    ([-\theta_0, \theta_0]\times [-\frac{1}{2}\log T, \infty))}
     \leqslant C < \infty,
    $$
where $C$ depends only on $\mu, k_i, \sigma$ and $\beta$ but not on $t, T$
and $u_0$.
\end{lemma}

The global existence of $w$ is not new, it has been obtained from the
existence of $r$ on $[0,T)$ in subsection 4.2. The estimate is
important and will be used below.


\subsection{Proof of Theorem \ref{thm:2}}
In this subsection we prove the existence and uniqueness of a discrete shrinking
self-similar solution on $[0,T)$.
Conditions \eqref{cond:k-shrink} and \eqref{cond:k-shrink-equi}
imply that $\tilde{g}_1 (s,w)$ and $\tilde{g}_2 (s,w)$ are $\log
b$-periodic in $s$. A similar discussion as in subsection 3.6 shows that
\eqref{prob-w} has a solution $\widetilde{P}(\theta,s)$, which is
$\log b$-periodic in $s$,
\begin{equation}\label{est of wtP}
\|\widetilde{P} (\theta ,s)\|_{C^{2+\mu, 1+ \mu/ 2}
([-\theta_0, \theta_0]\times \R)} \leqslant C (\mu,k_1, k_2, \sigma,\beta),
\end{equation}
and $\|w(\cdot, s) -\widetilde{P}(\cdot, s)\|_{C^2([-\theta_0, \theta_0])} \rightarrow
0$ as $s\rightarrow \infty$.

Now we recover $\widetilde{P}(\theta,s)$ back to a
corresponding solution of \eqref{prob}, that is, define $\widetilde{U}, \widetilde{\Xi}_i
\ (i=1,2)$ by
\begin{equation}\label{def:wtPwtXi}
\begin{array}{l}
e^{-s} e^{-\widetilde{P}} \cos\theta = \widetilde{U}(e^{-s} e^{-\widetilde{P}} \sin
\theta, T- e^{-2s}),\\
\widetilde{\Xi}_i (T- e^{-2s}) = (-1)^i e^{-s} e^{-\widetilde{P}(\theta_0, s)}
\sin\theta_0, \quad \mbox{ for } s\geqslant -\frac{1}{2}\log T.
\end{array}
\end{equation}
They are well-defined as in previous subsections. Moreover,
\begin{eqnarray*}
b^{-1} \widetilde{U}(e^{-s} e^{-\widetilde{P}(\theta,s)} \sin \theta, T- e^{-2s})
& = & b^{-1} e^{-s} e^{-\widetilde{P}(\theta, s)} \cos\theta  \\
&  = &  e^{-s -\log b} e^{-\widetilde{P}(\theta, s+\log b)}
\cos\theta \\
& = & \widetilde{U}(b^{-1} e^{-s} e^{-\widetilde{P}(\theta, s)} \sin \theta, T-
b^{-2} e^{-2s})
\end{eqnarray*}
and for $i=1,2$,
$$
b^{-1} \widetilde{\Xi}_i (T- e^{-2s}) = (-1)^i e^{-s-\log b} e^{-\widetilde{P}(\theta_0, s+\log b)}
\sin\theta_0  =  \widetilde{\Xi}_i (T- e^{-2s-2\log b}).
$$
Hence, for $t' := e^{2s} \in (0,T]$ and $-\widetilde{\Xi}_1 (T-t')\leqslant
x \leqslant \widetilde{\Xi}_2 (T-t')$ we have
\begin{equation}\label{exist-T}
b^{-1} \widetilde{U}(x,T-t') = \widetilde{U}(b^{-1}x, T-b^{-2} t'),\quad
b^{-1} \widetilde{\Xi}_i (T-t') = \widetilde{\Xi}_i (T-b^{-2} t')\ \ (i=1,2).
\end{equation}
This means that $(\widetilde{U},\widetilde{\Xi}_1, \widetilde{\Xi}_2)$ is a discrete shrinking
self-similar solution of \eqref{prob} on $[0,T)$.

We now prove the uniqueness result under the assumption that
$$
k_i (t,u) \equiv k_i (u)\ \quad (i=1,2).
$$

First, it is convenient to take a time shift and consider
the problem on $[-T,0)$. More precisely, as in section 1, we define
\begin{equation}\label{U-hat}
\wh{U}(x,t;T) := \widetilde{U}(x,T+t)\quad \mbox{for } -\wh{\Xi}_1
(t;T)\leqslant x \leqslant \wh{\Xi}_2(t;T),\ t\in [-T, 0),
\end{equation}
with
$$
\wh{\Xi}_i (t;T) := \widetilde{\Xi}_i (T+t)\quad \mbox{for }  t\in [-T,
0)\quad (i=1,2).
$$
Then \eqref{exist-T} implies that $\wh{U}$, $\wh{\Xi}_1$ and $\wh{\Xi}_2$
satisfy \eqref{self3} and \eqref{boundary3}. So $(\wh{U}, \wh{\Xi}_1, \wh{\Xi}_2)$
is a discrete self-similar solution of \eqref{prob-neg} on time-interval $t\in [-T,0)$,
and $\wh{U}, \wh{\Xi}_1, \wh{\Xi}_2$ all converge to $0$ as $t\to 0-0$.

Next, we construct self-similar solutions which decrease monotonically.
Indeed, as in subsection 3.8, under the assumption \eqref{cond:shrink}, we can choose
a concave initial data $u^*_0$ such that $a(u^*_{0x}) u^*_{0xx} <0$.
Let $u^*(x,t)$ be the solution of \eqref{prob} with initial data $u^*_0$,
then $u^*_t(x,t) <0$ by maximum principle, and $u^*(x,t)$ shrinks to
$0$ as $t\to T^*$ for some $T^*>0$. Moreover, for any given $T>0$, we can choose
$u^*_0$ sufficiently large such that $T^* >T$.

Converting this $u^*(x,t)$ to a new unknown $w^*(\theta,s)$ as in subsection 4.3
we have $1+ w^*_s (\theta,s)<0$ by \eqref{exp u w}, and so the $\omega$-limit $\widetilde{P}^*$ of
$w^*$ as in \eqref{est of wtP} satisfies $1+\widetilde{P}^*_s (\theta,s) \leqslant 0$.
Consequently, the corresponding functions $\widetilde{U}^*,\ \widetilde{\Xi}^*_1$
and $\widetilde{\Xi}^*_2$ defined by $\widetilde{P}^*$ as in \eqref{def:wtPwtXi} satisfy
$\widetilde{U}^*_t \leqslant 0$.  In addition, $\widetilde{U}^*$
is a shrinking self-similar solution of \eqref{prob} on $[0,T^*)$.
Hence,
\begin{equation}\label{U-hat-*}
\wh{U}^*(x,t;T^*) := \widetilde{U}^* (x,T^* +t),
\end{equation}
which is defined for $t\in [-T^*,0)$ and
$$
-\wh{\Xi}^*_1 (t) := -\widetilde{\Xi}_1 (T^* +t) \leqslant x\leqslant
\widetilde{\Xi}_2(T^* +t) =: \wh{\Xi}^*_2 (t),
$$
is a shrinking self-similar solution of \eqref{prob-neg} on $[-T^*,0)$.
By $\widetilde{U}^*_t \leqslant 0$ we have $\wh{U}^*_t \leqslant 0$.
By the strong maximum principle we even have
\begin{equation}\label{<0}
\wh{U}^*_t < 0,\quad \wh{\Xi}^*_{1t}<0, \quad \wh{\Xi}^*_{2t}<0 \ \quad \mbox{for}\ \ t\in [-T^*,0).
\end{equation}

\begin{lemma}\label{lem:equal}
Let $\wh{U}$ be as in \eqref{U-hat} and $\wh{U}^*$ be as in \eqref{U-hat-*}.
Then $\wh{U}(x,t) \equiv \wh{U}^*(x,t)$ for $t\in [-T,0)$.
\end{lemma}

\begin{proof}
Suppose on the contrary that, there exists $t_0 \in [-T,0)$ such that
\begin{equation}\label{contra 11}
\wh{U}(\cdot,t_0) \leqslant \wh{U}^* (\cdot,t_0) \quad \mbox{but} \quad \wh{U}(x,t_0)\not\equiv \wh{U}^*(x,t_0),
\end{equation}
or
\begin{equation}\label{contra 21}
\wh{U}^*(\cdot,t_0) \leqslant \wh{U}(\cdot,t_0) \quad \mbox{but} \quad \wh{U}^*(x,t_0)\not\equiv \wh{U}(x,t_0).
\end{equation}
Then we derive a contradiction from \eqref{contra 11} (In a similar way one can derive
a contradiction from \eqref{contra 21}).

By \eqref{contra 11} and by $\wh{U}^*_t <0$, there exists $\tau(t_0)\in [0, -t_0)$ such that
\begin{equation}\label{leq=}
\wh{U}(\cdot ,t_0) \preceq \wh{U}^* (\cdot,t_0 +\tau(t_0)).
\end{equation}
By comparison result Lemma \ref{comparison} we have
$$
\wh{U}\Big( \cdot, \frac{t_0}{b^2}\Big) = \wh{U}\Big( \cdot, t_0 + \frac{t_0}{b^2} -t_0\Big)
\ll \wh{U}^* \Big(\cdot, t_0 +\tau(t_0) + \frac{t_0}{b^2} -t_0 \Big) =
\wh{U}^* \Big(\cdot, \tau(t_0) + \frac{t_0}{b^2} \Big) .
$$
On the other hand, by \eqref{leq=} and the self-similarity \eqref{self3} we have
$$
\wh{U} \Big( \cdot, \frac{t_0}{b^2} \Big) \preceq \wh{U}^* \Big( \cdot, \frac{t_0 +\tau(t_0)}{b^2}\Big).
$$
Combining these inequalities with the fact that $\wh{U}^*_t <0$ we have
$\tau(t_0) > b^2 \tau(t_0)$, a contradiction. This proves the lemma.
\end{proof}

This completes the proof of Theorem \ref{thm:2}.

\subsection{Proof of Theorem \ref{thm:3}}
Under the assumption $k_i (t,u) \equiv k_i (u)\ (i=1,2)$, Theorem \ref{thm:3}
follows from the previous subsection easily.

Indeed, we can define positive functions $\Xi^\natural_1 (t)$ and $\Xi^\natural_2 (t)$
for $t\in (-\infty, 0)$ by
$$
\Xi^\natural_i (t) := \wh{\Xi}^*_i (t) \quad \mbox{for } t\in [-T^*,0), \ i=1,2,
$$
and define $U^\natural (x,t)$ on $(x,t)\in
[-\Xi^\natural_1 (t), \Xi^\natural_2(t)]\times (-\infty,0)$ by
$$
U^\natural (x,t) := \wh{U}^* (x,t) \quad \mbox{for }
x\in [-\Xi^\natural_1 (t),  \Xi^\natural_2(t)], \ t\in [-T^*,0).
$$
Lemma \ref{lem:equal} implies that the functions $U^\natural(x,t),\ \Xi^\natural_1 (t)$ and
$\Xi^\natural_2 (t)$ are well-defined (i.e., they do not depend on $T^*$).
The triple $(U^\natural, \;  \Xi^\natural_1,\; \Xi^\natural_2)$ is unique on time interval
$[-T,0)$ for any $T>0$, and so is unique in $(-\infty,0)$. Moreover, we have
$$
U^\natural_t (x,t) <0,\ \Xi^\natural_{1t} (t) <0,\ \Xi^\natural_{2t} (t) <0\quad \mbox{for } t<0
$$
by \eqref{<0}. This proves Theorem \ref{thm:3}.

\section*{Acknowledgments}
The author would like to thank Professors Hiroshi Matano and Ken-Ichi Nakamura
for helpful discussion. He also thanks the referees for valuable suggestions.

\medskip
Received xxxx 20xx; revised xxxx 20xx.
\medskip

\end{document}